\documentclass[12pt,reqno]{amsart}

\usepackage{amsmath,amsthm,amsfonts,amscd,amssymb,euscript,enumerate}
\numberwithin{equation}{section} 
\makeatletter
\@namedef{subjclassname@2020}{%
  \textup{2020} Mathematics Subject Classification}
\makeatother

\newcommand{\dbar}{\ensuremath{\bar \partial}}

\newcommand{\ad}{\ensuremath{\bar \partial^{*}  }}
\newcommand{\C}{\ensuremath{{\mathbb C}}}
\newcommand{\R}{\ensuremath{{\mathbb R}}}
\newcommand{\N}{\ensuremath{{\mathbb N}}}
\newcommand{\Q}{\ensuremath{{\mathbb Q}}}

\newcommand{\smooth}{\ensuremath{C^{\infty}}}
\newcommand{\smoothc}{\ensuremath{C_0 ^{\infty}}}
\newcommand{\ra}{\ensuremath{C^{\omega}}}
\newcommand{\I}{\ensuremath{\mathcal I}}
\newcommand{\II}{\ensuremath{\mathfrak I}}
\newcommand{\J}{\ensuremath{\mathcal J}}
\newcommand{\V}{\ensuremath{\mathcal V}}

\newcommand{\W}{\ensuremath{\mathcal W}}
\newcommand{\M}{\ensuremath{\mathcal M}}
\newcommand{\E}{\ensuremath{\mathcal E}}
\newcommand{\F}{\ensuremath{\mathcal F}}
\newcommand{\oka}{\ensuremath{\mathcal O}}
\newcommand{\nx}{\ensuremath{{\mathcal N}_x}}

\newcommand{\stb}{{\EuScript E}}

\newcommand{\An}{{\EuScript A}}

\newcommand{\gta}{{\mathfrak a}} 

\begin{document}
\title{The Kohn Algorithm on Denjoy-Carleman classes}
\author{Andreea C. Nicoara}

\address{School of Mathematics, Trinity College Dublin, Dublin 2, Ireland}

\email{anicoara@maths.tcd.ie}



\begin{abstract}
The equivalence of the Kohn finite ideal type and the D'Angelo finite type with the subellipticity of the $\dbar$-Neumann problem is extended to pseudoconvex domains in $\C^n$ whose defining function is in a Denjoy-Carleman quasianalytic class closed under differentiation. The proof involves algebraic geometry over a ring of germs of Denjoy-Carleman quasianalytic functions that is not known to be Noetherian and that is intermediate between the ring of germs of real analytic functions and the ring of germs of smooth functions. It is also shown that this type of ring of germs of Denjoy-Carleman functions satisfies the $\sqrt{acc}$ property, one of the strongest properties a non-Noetherian ring could possess.
\end{abstract}

\dedicatory{In memory of Professor Joseph J. Kohn}

\maketitle

\tableofcontents

\section{Introduction}

\bigskip The systematic study of the subellipticity of the
$\dbar$-Neumann problem on pseudoconvex domains in $\C^n$ was
initiated by Joseph J. Kohn in his paper \cite{kohnacta} published
in Acta Mathematica in 1979. Kohn defined subelliptic multipliers
for the $\dbar$-Neumann problem and showed they formed a
multiplier ideal sheaf. Kohn was the first to construct a
multiplier ideal sheaf, which has become since a standard object
in algebraic geometry. Kohn defined his subelliptic multipliers to
be germs of $\smooth$ functions and wrote down an algorithm that
generates an increasing chain of ideals of multipliers, whose
termination at the whole ring implies the subellipticity of the
$\dbar$-Neumann problem. This termination condition is called Kohn
finite ideal type. In the same paper \cite{kohnacta}, Joseph J.
Kohn proved a three-way equivalence for pseudoconvex domains in
$\C^n$ with real analytic $\ra$ boundary, namely that the
subellipticity of the $\dbar$-Neumann problem on $(p,q)$ forms is
equivalent to Kohn finite ideal type for $(p,q)$ forms and is also
equivalent to the condition that holomorphic varieties of complex
dimension $q$ have finite order of contact with the boundary of
the domain. The latter condition is called finite D'Angelo type
after John D'Angelo who investigated its properties in detail in \cite{opendangelo} and proved the crucial fact that it is an open condition.

The next important development in the investigation of the
subellipticity of the $\dbar$-Neumann problem came in the mid
1980's with a series of three deep papers by David Catlin,
\cite{catlinnec}, \cite{catlinbdry}, and \cite{catlinsubell}, in
which he proved the equivalence of two out of the three properties
that appear in Kohn's theorem for real analytic domains, namely
that for a smooth, pseudoconvex domain in $\C^n$ the
subellipticity of the $\dbar$-Neumann problem is equivalent to
another notion of finite type that is effectively related to D'Angelo finite type; see \cite{bazilandreea2019}. Catlin's construction involves a very careful analysis of the vanishing orders of the defining function of the domain in different directions but not
address the question of whether the Kohn algorithm terminates, i.e. whether Kohn finite
ideal type is also equivalent to the other two properties for
smooth, pseudoconvex domains in $\C^n.$ This three-way equivalence
was subsequently posited by Joseph J. Kohn and David Catlin and is called the Kohn
conjecture.

The current investigation of the equivalence of types for
Denjoy-Carleman quasianalytic classes is an intermediate case between the $\ra$ one settled by Kohn in 1979 and the Kohn conjecture. It involves classes of functions without necessarily convergent Taylor expansions that still satisfy the {\L}ojasiewicz
inequalities, properties crucial to the
equivalence of types. Denjoy-Carleman quasianalytic classes are not standardly used in several complex variables. Besides generalizing Kohn's 1979 result, the emphasis here is to highlight how far Kohn's ideas from the real-analytic case can be stretched and what properties would be necessary to prove the Kohn conjecture. More recently, the Kohn algorithm in the context of the Denjoy-Carleman quasianalytic classes was treated in \cite{LeviCore2} without giving a proof of termination, however.

The Denjoy-Carleman quasianalytic classes are subrings of
the ring of smooth functions such that certain bounds hold on their
derivatives of various orders. The sequence of bounds behaves
according to the Denjoy-Carleman Theorem, which means that the
Taylor morphism is injective on any such Denjoy-Carleman
quasianalytic class, preventing flat functions from lying therein. Each Denjoy-Carleman quasianalytic class considered here will be strictly larger than the class of real analytic functions, thus also containing functions with non-convergent Taylor expansions. Denjoy-Carleman quasianalytic classes do not satisfy
the Weierstrass Division Property as Childress proved in
\cite{childress} or the Weierstrass Preparation Theorem as proven in \cite{nonwpt}, and yet Edward Bierstone and Pierre Milman established in \cite{bm} that resolution of singularities holds. As a result of the Bierstone-Milman construction, all three {\L}ojasiewicz
inequalities follow along with topological Noetherianity, a weaker
condition than Noetherianity but of great use for the equivalence
of types considered here. \cite{bm} assumes sheaves are of finite type. More recently, this hypothesis was relaxed to quasicoherence in \cite{bmv}.

We shall show in this paper that
for a pseudoconvex domain in $\C^n$ with boundary in any such
Denjoy-Carleman quasianalytic class, the
three-way equivalence proved by Kohn for pseudoconvex $\ra$
domains extends to this case:

\medskip
\newtheorem{mainthm}{Main Theorem}[section]
\begin{mainthm}
Let $\Omega$ in $\C^n$ be a pseudoconvex domain with boundary in a
Denjoy-Carleman quasianalytic class $C_M (\overline{\Omega})$
defined in a neighborhood of $\overline{\Omega},$ where the
sequence of bounds $M$ is such that $C_M (\overline{\Omega})$ is
closed under differentiation. Let $x_0 \in b \Omega$ be any point
on the boundary of the domain, and let $\tilde U$ be an
appropriately small neighborhood around $x_0$. The following three
properties are equivalent: \label{maintheorem}
\begin{enumerate}
\item[(i)] The $\dbar$-Neumann problem for $(p,q)$ forms is subelliptic on $\tilde U$;
\item[(ii)] The Kohn algorithm on $(p,q)$ forms terminates at $x_0$ by generating the entire ring of germs of smooth functions on $b \Omega$ at $x_0,$ $\smooth_{b \Omega} (x_0)$;
\item[(iii)] The order of contact of holomorphic varieties of complex dimension $q$ with
the boundary of the domain $\Omega$ in $\tilde U$ is finite.
\end{enumerate}
\end{mainthm}

As mentioned at the beginning of this introduction, the
implication (ii) $\implies$ (i) was already done by Joseph J. Kohn
for $\smooth$ functions in \cite{kohnacta}, so his construction
applies here as well. The implication (i) $\implies$ (iii) is the
contrapositive of Catlin's Theorem 1 in \cite{catlinnec}. We
should add that a complex variety passing through $x_0$ that has
infinite order of contact with $b \Omega$ is already sitting
inside $b \Omega$ since the Denjoy-Carleman quasianalytic
functions considered here contain no flat functions. Thus, the aim
of this paper is to prove the implication (iii) $\implies$ (ii).
The proof that will be given is purely qualitative in nature.
There is no attempt to compute an effective bound for the
subelliptic gain in the $\dbar$-Neumann problem in terms of the
D'Angelo type and the dimension. The focus is instead on the
interaction of the algebraic properties of the Denjoy-Carleman
quasianalytic classes with the algebra of the Kohn algorithm. The most delicate point is showing that the sheaves involved in the argument are quasicoherent, which is achieved by combining results and techniques from \cite{abn} and \cite{andreeaqf}. For results and discussion of effective results in the sum of squares of holomorphic case and the real analytic case, the reader should consult \cite{catlindangelononeff}, \cite{dangelo95}, \cite{kimzaitsev2}, \cite{kimzaitsev}, \cite{racase}, \cite{siu2010}, and \cite{siu2017}. In fact, Theorem~\ref{maintheorem} is true in an even more general setting than the Denjoy-Carleman quasianalytic classes. For any class of functions satisfying axioms (3.1)-(3.6) on p.3-4 of \cite{bm}, the resolution of singularities goes through and with it all algebraic geometric properties necessary for Theorem~\ref{maintheorem}. An example of more general quasianalytic classes with such properties are classes of $\smooth$ functions that are definable in a given polynomially bounded o-minimal structure. The reader is directed to \cite{rsw} for more details on these classes.

It is widely believed that the rings of germs of Denjoy-Carleman classes are non-Noetherian. Nonetheless, the topological Noetherianity proven by Bierstone and Milman permits the construction of certain finitely generated subideals of the ideals of multipliers via which it can shown that the ascending chain in the Kohn algorithm stabilizes. Likewise, the presence of {\L}ojasiewicz
inequalities enables a Nullstellensatz for a notion of real radical particular to the Kohn algorithm. Therefore, the Denjoy-Carleman setting provides a glimpse into the algebra necessary for the direct
proof of the Kohn conjecture.

In addition to the main result, Theorem \ref{maintheorem}, concerning the equivalence of types, we derive a very simple consequence of Bierstone's and Milman's work in \cite{bm}, namely that all the rings of germs corresponding to Denjoy-Carleman classes closed under differentiation satisfy the $\sqrt{acc}$ condition for the real radical notion $\sqrt[\R]{\,}.$ In commutative algebra, this condition is the nicest a non-Noetherian ring can satisfy. The reader may consult \cite{kaplansky} and \cite{op} for additional information on this type of rings.

\smallskip
\newtheorem{radacc}[mainthm]{Theorem}
\begin{radacc}
\label{radaccprop} Let $x_0$ be a point in $\R^n$ or $\C^n,$ and let $C_M (x_0)$ be a ring of germs of Denjoy-Carleman quasianalytic functions closed under differentiation. The local ring $C_M (x_0)$ has the $\sqrt{acc}$ property. In other words, if $\I_1 \subset \I_2 \subset \cdots$ is an ascending chain of ideals in $C_M (x_0),$ then the ascending chain of radical ideals $\sqrt[\R]{\I_1} \subset \sqrt[\R]{\I_2} \subset \cdots$ stabilizes, i.e. there exists a $k$ such that $\sqrt[\R]{\I_j}=\sqrt[\R]{\I_k}$ for all $j \geq k.$
\end{radacc}

The paper is organized as follows: Section~\ref{DCdescription}
introduces the Denjoy-Carleman quasianalytic classes and outlines
a number of their properties. At the end of this section, Theorem~\ref{radaccprop} is proven. Section~\ref{subkohnalg} recalls the
Kohn algorithm and other matters related to the subellipticity of
the $\dbar$-Neumann problem. Finally, the equivalence of types,
Theorem~\ref{maintheorem}, is proven in Section~\ref{equivtypepf}.

\bigskip\bigskip

\section{The Denjoy-Carleman Quasianalytic Classes}
\label{DCdescription}

\bigskip
The set-up described here comes from the Bierstone and Milman paper
\cite{bm}. Another excellent source for properties of more general quasianalytic local rings is Vincent Thilliez's
expository paper \cite{thilliez}.

We start with the definition of a quasianalytic class. Let us make
the identification $\C^n \simeq \R^{2n},$ where $z=(z_1, \dots,
z_n)=(x_1+i x_{n+1},\dots, x_n+i x_{2n}).$

\smallskip
\newtheorem{quasian}{Definition}[section]
\begin{quasian}
Let $U$ be a connected open set in $\C^n$, and let $M=\{M_0, M_1,
M_2, \dots \}$ be an increasing sequence of positive real numbers,
where $M_0=1.$ $C^\R_M (U)$ consists of all $\R$-valued $f \in
\smooth (U)$ satisfying that for every compact set $K \subset U$,
there exist constants $A, B>0$ such that
\begin{equation}
\left|\frac{1}{\alpha !}\, D^\alpha  f(x)\right| \leq A \,
B^{|\alpha|} \, M_{|\alpha|} \label{realquasianest}
\end{equation}
for any $x \in K,$ where $\alpha$ is a multi-index in $\N^{2n}$,
$D^\alpha = \frac{\partial^{\, |\alpha|}}{\partial x_1^{\,
\alpha_1} \cdots \,
\partial x_{2n}^{\, \alpha_{2n}}}.$ $C^\R_M (U)$ is called
quasianalytic if the Taylor morphism assigning to each $f \in
C^\R_M (U)$ its Taylor expansion at $a \in U$ is injective for all
$a \in U.$ \label{realquasiandef}
\end{quasian}

\smallskip\noindent Since the $\dbar$-Neumann problem is posed on
$(p,q)$ forms in $\C^n,$ we must also introduce another class $C_M
(U)$ of $\C$-valued quasianalytic functions. We will examine the properties
of both classes $C^\R_M (U)$ and $C_M (U),$ going back and forth
between them as necessary.

\smallskip
\newtheorem{Cquasian}[quasian]{Definition}
\begin{Cquasian}
Let $U$ be a connected open set in $\C^n$, and let $M=\{M_0, M_1,
M_2, \dots \}$ be an increasing sequence of positive real numbers,
where $M_0=1.$ $C_M (U)$ consists of all $\C$-valued $f \in
\smooth (U)$ such that $f=g+ih$ for $g$ and $h$ real valued and $g, h \in C^\R_M (U).$ $C_M (U)$ is called
quasianalytic if the Taylor morphism assigning to each $f \in C_M
(U)$ its Taylor expansion at $a \in U$ is injective for all $a \in
U.$ \label{Cquasiandef}
\end{Cquasian}

\smallskip\noindent {\bf Remarks:}

\smallskip\noindent (1) Clearly, $$f=g+ih \in C_M (U) \quad\iff\quad
g,h \in C^\R_M(U)$$ and $$C_M(U) \quad \text{quasianalytic} \quad\iff\quad
C^\R_M(U) \quad \text{quasianalytic}.$$

\smallskip\noindent (2) Consider any open set $\tilde U \subset U.$
If $f \in C_M (U),$ then it is obvious from
Definition~\ref{Cquasiandef} that the restriction $f\big|_{\tilde
U} \in C_M (\tilde U).$ The same holds for $C^\R_M (U).$ This observation naturally leads to the definition of germs of Denjoy-Carleman functions.

\smallskip
\newtheorem{Cquasiangerms}[quasian]{Definition}
\begin{Cquasiangerms}
Let $x_0$ be a point in $\C^n$. The set of germs $C_M(x_0)$ (respectively $C^\R_M (x_0)$) of Denjoy-Carleman functions at $x_0$  consists of all germs of $\C$-valued (respectively $\R$-valued) functions $f \in
C_M (U)$ (respectively $f \in C^\R_M (U)$) for $U$ a neighborhood of $x_0.$ $C_M (x_0)$ (respectively $C^\R_M (x_0)$)  is called quasianalytic if the Taylor morphism assigning to each $f \in C_M
(x_0)$ (respectively $f \in C^\R_M (x_0)$) its Taylor expansion at $x_0$ is injective. 
\label{Cquasiangermsdef}
\end{Cquasiangerms}

\bigskip\noindent In order to ensure closure under multiplication that makes $C^\R_M(U),$ $C_M(U),$ $C^\R_M (x_0),$ and $C_M (x_0)$ into rings, which is essential for the algebraic considerations that
will follow, we must impose an additional condition on the sequence
$M$ of bounds:

\smallskip
\newtheorem{logconvex}[quasian]{Definition}
\begin{logconvex}
The sequence $M=\{M_0, M_1, M_2, \dots \}$ is called
logarithmically convex if
\begin{equation}
\frac{M_{j+1}}{M_j}\leq\frac{M_{j+2}}{M_{j+1}}
\label{logcon}
\end{equation}
for all $j \geq
0,$ i.e. the sequence of subsequent quotients is increasing.
\label{logconvexdef}
\end{logconvex}

\smallskip\noindent Definition~\ref{logconvexdef} and the assumption $M_0
=1$ together imply that $\{\, (M_j)^{\, \frac{1}{j}}\, \}_{j \geq 1}$ is an
increasing sequence of positive real numbers greater than or equal
to $1.$ Real analytic functions satisfy estimate
~\eqref{realquasianest} with $M_{|\alpha|}=1$ for all $\alpha.$ Therefore, all $\R$-valued real analytic
functions belong to $C^\R_M (U)$ and all $\C$-valued real analytic
functions belong to $C_M (U).$ Fortunately, there exists a condition on $M$ in the classical literature that precisely identifies when $C^\R_M (U)$ equals the ring of $\R$-valued real analytic functions on $U,$
$\ra (U);$  see Corollary 1 of Theorem 1 in \cite{thilliez}.

\smallskip
\newtheorem{equaltora}[quasian]{Proposition}
\begin{equaltora}
$C^\R_M (U) \, = \, \ra (U)$ iff $\sup_{j \geq 1} (M_j)^{\,
\frac{1}{j}} < \infty.$
\end{equaltora}

\smallskip\noindent To ensure both $C^\R_M (U)$ and $C_M (U)$ are
strictly larger than the rings of $\R$-valued real analytic
functions and $\C$-valued real analytic functions respectively, we shall assume
\begin{equation}
\lim_{j \rightarrow \infty} (M_j)^{\, \frac{1}{j}} = \infty.
\label{nonra}
\end{equation}
This strict inclusion property then propagates to the rings of germs $C^\R_M (x_0)$ and $C_M (x_0).$

 \'{E}mile Borel first constructed quasianalytic functions that were not also real analytic in \cite{borel1} and \cite{borel2}. Other examples can be found in \cite{thilliez}. Holmgren's work on the heat equation subsequently showed that solutions to certain partial differential equations belong to classes of functions intermediate between real analytic and $\smooth.$ In response to Holmgren's work, Jacques Hadamard then asked in 1912 whether there exists a condition on the sequence $M$ that forces the corresponding class to be quasianalytic. The Denjoy-Carleman Theorem provides the answer:

\smallskip
\newtheorem{dcthm}[quasian]{Denjoy-Carleman Theorem}
\begin{dcthm}
Let the sequence $M=\{M_0, M_1, M_2, \dots \}$ be logarithmically
convex, then $C^\R_M(U),$ $C_M(U),$ $C^\R_M (x_0),$ and $C_M (x_0)$ are quasianalytic iff
$\:\sum_{k=0}^\infty \: \frac{M_k}{(k+1) \, M_{k+1}} = \infty.$ \label{dctheorem}
\end{dcthm}

\smallskip\noindent Such a class $C^\R_M(U),$ $C_M(U),$ $C^\R_M (x_0),$ or $C_M (x_0)$ satisfying the Denjoy-Carleman Theorem is called a Denjoy-Carleman quasianalytic
class. We have just one more condition to impose on $M,$ and we
will have described completely the Denjoy-Carleman quasianalytic
classes on which we will consider the equivalence of types. This
condition due to Szolem Mandelbroit guarantees that $C^\R_M(U),$ $C_M(U),$ $C^\R_M (x_0),$ and $C_M (x_0)$ are all
closed under differentiation:

\smallskip
\newtheorem{derivok}[quasian]{Proposition}
\begin{derivok}
$C^\R_M(U),$ $C_M(U),$ $C^\R_M (x_0),$ and $C_M (x_0)$ are closed under differentiation
iff
\begin{equation}
\sup_{j \geq 1} \left( \frac{M_{j+1}}{M_j}\right)^{\, \frac{1}{j}}
< \infty. \label{derivokcond}
\end{equation}
\end{derivok}

\smallskip\noindent {\bf Remark:} Clearly, $C^\R_M (U)$ closed under differentiation implies $C_M (U)$ is also closed under differentiation and the same for $C^\R_M (x_0)$ and $C_M (x_0).$

\smallskip\noindent We shall summarize now the results of Section
4 of \cite{bm} in a proposition that lists all the properties of a
Denjoy-Carleman quasianalytic class of the type described above
that enable the resolution of singularities. This result is only
stated for $C^\R_M (U)$ because some of the properties contained
therein are awkward to give in complex coordinates. As mentioned in the introduction, whenever these properties hold, the resolution of singularities goes through, so technically in a more general setting than we are using here.

\smallskip
\newtheorem{allthere}[quasian]{Proposition}
\begin{allthere}
Let $C^\R_M (U)$ be a Denjoy-Carleman quasianalytic class that
also satisfies conditions ~\eqref{logcon}, ~\eqref{nonra}, and ~\eqref{derivokcond}
on its sequence of bounds $M.$ $C^\R_M (U)$ has the following
properties: \label{allpropsthere}
\begin{enumerate}
\item[(i)] $C^\R_M (U)$ contains all real analytic functions on $U$
as well as the restrictions of all polynomials on $U,$ ${\mathcal
P}(U);$
\item[(ii)] $C^\R_M (U)$ is closed under composition, namely if $U$
and $V$ are open subsets of $\C^n$ and $\C^p$ respectively, $f \in
C^\R_M (V),$ and $g = (g_1, \dots, g_{2p}): U \longrightarrow V$
is such that $g_j \in C^\R_M (U)$ for $1 \leq j \leq 2p,$ then $f
\circ g \in C^\R_M (U);$
\item[(iii)] $C^\R_M (U)$ is closed under differentiation;
\item[(iv)] $C^\R_M (U)$ is closed under division by a coordinate,
i.e. if $f \in C^\R_M (U)$ and $$f(x_1, \dots, x_{i-1}, a_i,
x_{i+1}, \dots, x_{2n}) \equiv 0,$$ then there exists $h \in
C^\R_M (U)$ such that $f(x)= (x_i - a_i) \, h(x);$
\item[(v)] $C^\R_M (U)$ is closed under inverse, namely let $U$ and
$V$ be open subsets of $\C^n$ and let $\varphi = (\varphi_1,
\dots, \varphi_{2n}): U \longrightarrow V$ be such that $\varphi_i
\in C^\R_M (U)$ for $1 \leq i \leq 2n, $ $a \in U,$
$\varphi(a)=b,$ and the Jacobian matrix $$\frac{\partial
\varphi}{\partial x}(a)=\frac{\partial (\varphi_1, \dots,
\varphi_{2n}) }{\partial (x_1,\dots, x_{2n})}(a)$$ be invertible,
then there exist neighborhoods $U'$ of $a$ and $V'$ of $b$ as well
as a mapping $\psi = (\psi_1, \dots, \psi_{2n}): V'
\longrightarrow U'$ such that $\psi_i \in C^\R_M (V')$ for $1 \leq
i \leq 2n,$ $\psi(b)=a$ and $\varphi \circ \psi$ is the identity
mapping on $V'.$
\end{enumerate}
\end{allthere}

\smallskip\noindent {\bf Remarks:}

\noindent (1) Proposition~\ref{allpropsthere} part (v) is
equivalent to the Implicit Function Theorem, which can be stated
in this context as follows: Let $U$ be an open subset of $\C^n
\times \C^p$ with product coordinates $(x,y)=(x_1, \dots, x_{2n},
y_1, \dots, y_{2p}).$ Let $f_1, \dots, f_{2p} \in C^\R_M (U),$
$(a,b) \in U,$ $f(a,b)=0=(f_1(a,b), \dots, f_{2p}(a,b)),$ and
$\frac{\partial f}{\partial y} (a,b)$ be invertible, then there
exists a product neighborhood $V \times W$ of $(a,b)$ in $U$ and a
mapping $g=(g_1, \dots, g_{2p}): V \longrightarrow W$ such that
$g_i \in C^\R_M (V)$ for all $1 \leq i \leq 2p,$ $g(a)=b,$ and
$f(x,g(x))=0$ for all $x \in V.$

\noindent (2) Proposition~\ref{allpropsthere} part (v) implies
that $C^\R_M (U)$ is closed under reciprocal, i.e. if $f \in
C^\R_M (U)$ is such that $f(x) \neq 0$ for all $x \in U,$ then
$\frac{1}{f} \in C^\R_M (U).$

\noindent (3) Parts (ii) and (iv) of Proposition~\ref{allpropsthere} imply part (iii); see \cite{bmv}.

\noindent (4) Using the properties of $C^\R_M (U)$ in
Proposition~\ref{allpropsthere}, we can construct manifolds of
class $C^\R_M$ as well as functions on such manifolds. Since
real and imaginary parts of functions in $C_M (U)$ are in $C^\R_M
(U),$ we automatically also get well-defined manifolds of class
$C_M$ as well as functions on those manifolds.

\bigskip As it will be seen in Section~\ref{subkohnalg},
the Kohn algorithm yields an increasing chain of ideals of germs of
functions defined on the neighborhood of a point in the boundary of the domain $\Omega \subset\C^n.$ We thus have to explore here what it means to have germs of
functions in the class $C_M$ defined on a manifold of the same
class as well as sheaves of ideals of functions in the same class. Let
$X$ be a manifold of class $C_M$ such that $X \subset U
\subset \C^n.$ Since $X$ is embedded in $\C^n,$ we equip it with
the subset topology, where $\C^n$ has the standard Euclidean
metric topology on it. Let ${\mathfrak G}_X^{C_M}$ be the sheaf
of germs of functions in the class $C_M$ defined at points of $X.$ We will denote the stalk of this sheaf at a point $x_0 \in X$ as $C_{M,X} (x_0).$ We now
consider a sheaf of ideals $\II \subset {\mathfrak G}_X^{C_M}.$ The variety corresponding to this sheaf is formally defined as $$\V(\II)= supp \: \frac{ {\mathfrak G}_X^{C_M}}{\II}.$$ The reader should note that $\V(\II)$ consists of all points $a \in X,$ where the ideal $\II_a$ does not contain any unit in the ring of germs of functions of class $C_M$ on $X$ at $a,$ which we denoted $C_{M,X} (a).$ This definition is consistent with that in \cite{bm}. At times, we will need to work with sheaves of ideals of finite type. We give the same definition here as in \cite{bm}.

\smallskip
\newtheorem{finitetype}[quasian]{Definition}
\begin{finitetype}
A sheaf of ideals $\II \subset {\mathfrak G}_X^{C_M}$ is said
to be of finite type if for each $a \in X,$ there exist a
neighborhood $U_a$ of $a$ in $X$ and finitely many sections $f_1,
\dots, f_p \in C_M (U_a)$ such that for all $b \in U_a$ the stalk
$\II_b$ is generated by the germs of $f_1, \dots, f_p$ at $b.$
\end{finitetype}

\smallskip\noindent For a manifold $X$ of class $C^\R_M $ with germs
${\mathfrak G}_X^{C^\R_M} $ on it, the definition of a variety corresponding to a sheaf of ideals and that of a sheaf of ideals of finite type are the same as in the case of $C_M.$ It shall also be noted here that an ideal of finite type in a ring means a finitely generated ideal.

We also need to introduce here the notion of the germ of a variety corresponding to an ideal of germs of Denjoy-Carleman functions that is of finite type. This will allow us to visualize the varieties obtained pointwise at the various steps of the Kohn algorithm. We are assuming the rings of germs $C^\R_M (x_0),$ $C_M (x_0),$ and $C_{M,X} (x_0)$ all satisfy the Denjoy-Carleman Theorem, Theorem~\ref{dctheorem}. This means that if two functions $f$ and $g$ defined on a neighborhood of $x_0$ coincide on a perhaps smaller neighborhood of $x_0,$ then they must be the same function as they would have the same Taylor expansion at $x_0$ itself. Therefore, a germ of a Denjoy-Carleman function is an equivalence class containing one and only one element, so whenever we are given a germ, we can work with its representative in a neighborhood of $x_0.$ This is the same as in the real analytic case. Consider now an ideal $\I$ of finite type in any of the rings of germs $C^\R_M (x_0),$ $C_M (x_0),$ or $C_{M,X} (x_0),$ and let $\I=(f_1, \dots, f_p).$ Each generator $f_j$ is a germ corresponding to a Denjoy-Carleman function, which we will also denote by $f_j$ since it is unique,  and $f_j$ is defined on some neighborhood of $x_0,$ which we will denote by $U_j.$ Consider now the intersection of neighborhoods $U= U_1 \cap \dots \cap U_p.$ All $f_j$'s are defined on $U,$ so the germ of the variety $V(\I)$ corresponding to $\I$ is the set of common zeroes of $f_1, \dots, f_p$ in $U.$ Certain arguments will require us to shrink this $U$ as it is standard when working with germs.

For more sophisticated arguments, we will need to consider quasicoherent sheaves of ideals.

\smallskip
\newtheorem{quasicoh}[quasian]{Definition}
\begin{quasicoh}
\label{quasicoherent} A sheaf of ideals $\II \subset {\mathfrak G}_X^{C_M}$ is said
to be quasicoherent if for each $a \in X,$ there exist a
neighborhood $U_a$ of $a$ in $X$ and a family of sections $\{f_\alpha\} \in C_M (U_a)$ not necessarily finite such that for all $b \in U_a$ the stalk
$\II_b$ is generated by the germs of $\{f_\alpha\}$ at $b.$
\end{quasicoh}

\smallskip We can now state two of the consequences of applying
resolution of singularities to the Denjoy-Carleman quasianalytic
class $C^\R_M ,$ which Edward Bierstone and Pierre Milman
obtained in \cite{bm}. We shall state as corollaries in both
instances the corresponding results for $C_M,$ which are the
ones we are more interested in here. Both of these are crucial for the
equivalence of types on domains of class $C_M.$ The first of
these important consequences is topological Noetherianity, which
is Theorem 6.1 of \cite{bm} for sheaves of ideals of finite type and also follows for quasicoherent sheaves of ideals from the resolution of singularities argument in \cite{bmv}.

\smallskip
\newtheorem{rtopnoet}[quasian]{Theorem (Topological Noetherianity)}
\begin{rtopnoet}
Let $\I_1, \I_2, \dots$ be any sequence of quasicoherent sheaves of ideals in ${\mathfrak G}_X^{C^\R_M }$ such that the corresponding
varieties form a decreasing sequence $\V (\I_1) \supseteq \V
(\I_2) \supseteq \cdots.$ Given some compact set $K$ in the
topology of $X,$ the sequence $\V (\I_1) \supseteq \V (\I_2)
\supseteq \cdots$ stabilizes in a neighborhood of $K,$ i.e. there
exists some $k$ such that in a neighborhood of $K,$ $\V (\I_j) =
\V (\I_k)$ for all $j \geq k.$ \label{rtopnoetherianity}
\end{rtopnoet}

\smallskip\noindent After trivial modifications to the proof, Theorem 6.1 of \cite{bm} of Bierstone and Milman also holds for germs as Federica Pieroni observed in \cite{pieroni}:


\smallskip
\newtheorem{grtopnoet}[quasian]{Theorem (Topological Noetherianity for Germs of Class $C^\R_M$)}
\begin{grtopnoet}
Let $\I_1, \I_2, \dots$ be any sequence of ideals of finite type
in $C^\R_M (x_0)$ for $x_0 \in \R^n$ or $x_0 \in \C^n$ or in $C^\R_{M,X}(x_0)$ for $X$ a manifold of class $C^\R_M$ and $x_0 \in X$ such that the corresponding germs of 
varieties form a decreasing sequence $\V (\I_1) \supseteq \V
(\I_2) \supseteq \cdots.$ The sequence $\V (\I_1) \supseteq \V (\I_2)
\supseteq \cdots$ stabilizes, i.e. there
exists some $k$ such that $\V (\I_j) =
\V (\I_k)$ for all $j \geq k.$ \label{grtopnoetherianity}
\end{grtopnoet}

\smallskip\noindent To any ideal $\I \subset C_M(x_0)$ can be associated an ideal $\I^\R \subset C^\R_M(x_0)$ defined by $$\I^\R = \{f \in  C^\R_M (x_0)\:\:
| \:\: \exists \, g \in \I \:\text{s.t.}
\:f=Re \, g \: \text{or} \: f=Im \, g \}.$$ Clearly, $\I^\R$ is of finite type if $\I$ is and $\V(\I^\R)=\V(\I).$ The same construction applies to some ideal $I \subset C_{M,X}(x_0)$  for $X$ a manifold of class $C_M$ and $x_0 \in X.$ We thus obtain:

\smallskip
\newtheorem{topnoet}[quasian]{Corollary}
\begin{topnoet}
Let $\I_1, \I_2, \dots$ be any sequence of ideals of finite type
in $C_M (x_0)$ for $x_0 \in \R^n$ or $x_0 \in \C^n$ or in $C_{M,X}(x_0)$  for $X$ a manifold of class $C_M$ and $x_0 \in X$ such that the corresponding germs of 
varieties form a decreasing sequence $\V (\I_1) \supseteq \V
(\I_2) \supseteq \cdots.$ The sequence $\V (\I_1) \supseteq \V (\I_2)
\supseteq \cdots$ stabilizes, i.e. there
exists some $k$ such that $\V (\I_j) =
\V (\I_k)$ for all $j \geq k.$ \label{topnoetherianity}
\end{topnoet}

\smallskip\noindent The second of the important consequences of
Edward Bierstone's and Pierre Milman's work is the full set of
{\L}ojasiewicz inequalities for the Denjoy-Carleman quasianalytic
class $C^\R_M,$ Theorem 6.3 of \cite{bm}.

\smallskip
\newtheorem{rloj}[quasian]{Theorem}
\begin{rloj}
The three {\L}ojasiewicz inequalities hold for functions of class $C^\R_M $ as
follows:
\begin{enumerate}
\item[(I)] Let $X$ be a manifold of class $C^\R_M,$ and let
$f, g \in C^\R_M (X),$ i.e. these functions are defined in a
neighborhood of $X.$ If $\{ x \in X \: \big| \: g(x)=0 \}
\subseteq \{ x \in X \: \big| \: f(x)=0 \}$ in a neighborhood of a
set $K$ compact in the topology of $X,$ then there exist $C,
\alpha >0$ such that
$$|g(x)| \geq C \, |f(x)|^{\, \alpha}$$ for all $x$ in a neighborhood
of $K.$ Furthermore, $\inf \alpha$ is a positive rational number.
\item[(II)] Let $f \in C^\R_M (U)$ for $U$ an open set such that $U \subset \R^n,$ and set $Z = \{ x \in U \: \big|
\: f(x)=0 \}.$ Suppose that $K \subset U$ is compact. Then there
exist $C >0$ and $\nu \geq 1$ such that $$|f(x)| \geq C \,
d(x,Z)^{\, \nu}$$ in a neighborhood of $K,$ where $d(\, \cdot \, ,
Z)$ is the Euclidean distance to $Z$ and $\inf \nu \in \Q.$
\item[(III)] Let $f \in C^\R_M (U),$ and let $K$ be a compact subset
of $U$ such that $\nabla_\R \, f(x) =0$ only if $f(x) =0,$ where
$\nabla_\R$ is the gradient in real coordinates $(x_1, \dots,
x_{2n}).$ Then there exist $C
>0$ and $\mu$ satisfying $0< \mu \leq 1$ such that $$|\nabla_\R \,
f(x)| \geq C \, |f(x)|^{\, 1 - \mu}$$ in a neighborhood of $K$ and
$\sup \mu \in \Q.$
\end{enumerate} \label{rlojasiewicz}
\end{rloj}

\smallskip\noindent Since the real and imaginary parts of functions in
$C_M (U)$ are in $C^\R_M (U),$ it is clear that parts (I) and (II)
will also be true for the class $C_M.$ Part (III) is irrelevant for type
equivalence, so we will not restate it for $C_M.$

\smallskip
\newtheorem{loj}[quasian]{Corollary}
\begin{loj}
The first two of the {\L}ojasiewicz inequalities above hold for functions of class
$C_M$:
\begin{enumerate}
\item[(I)] Let $X$ be a manifold of class $C_M,$ and let
$f, g \in C_M (X),$ i.e. these functions are defined in a
neighborhood of $X.$ If $\{ x \in X \: \big| \: g(x)=0 \}
\subseteq \{ x \in X \: \big| \: f(x)=0 \}$ in a neighborhood of a
set $K$ compact in the topology of $X,$ then there exist $C,
\alpha >0$ such that
$$|g(x)| \geq C \, |f(x)|^{\, \alpha}$$ for all $x$ in a neighborhood
of $K.$ Furthermore, $\inf \alpha$ is a positive rational number.
\item[(II)] Let $f \in C_M (U)$ for $U$ an open set such that $U \subset \R^n,$ and set $Z = \{ x \in U \: \big|
\: f(x)=0 \}.$ Suppose that $K \subset U$ is compact. Then there
exist $C >0$ and $\nu \geq 1$ such that $$|f(x)| \geq C \,
d(x,Z)^{\, \nu}$$ in a neighborhood of $K,$ where $d(\, \cdot \, ,
Z)$ is the Euclidean distance to $Z$ and $\inf \nu \in \Q.$
\end{enumerate} \label{lojasiewicz}
\end{loj}

\smallskip\noindent {\bf Proof:} Note that $f=f_1+if_2$ and $g=g_1+ig_2$ for $f_1, f_2, g_1, g_2 \in C^\R_M (X).$ To prove part (I), we look at $g_1^2+g_2^2$ and compare it to $f_1^2+f_2^2.$ Clearly, $\{ x \in X \: \big| \: f(x)=0 \} = \{ x \in X \: \big|\: f^2_1+f^2_2=0 \}$ and $\{ x \in X \: \big| \: g(x)=0 \} = \{ x \in X \: \big|\: g^2_1+g^2_2=0 \}$ so $ \{ x \in X \: \big|\: g^2_1+g^2_2=0 \} \subseteq  \{ x \in X \: \big|\: f^2_1+f^2_2=0 \}.$ The result follows from part (I) of Theorem~\ref{rlojasiewicz}. Similarly, part (II) follows from part (II) of Theorem~\ref{rlojasiewicz} applied to $f_1^2+f_2^2,$ which is $\R$-valued and has the same zero set as $f.$ For part (II), $f_1, f_2 \in C^\R_M(U).$  \qed

\smallskip\noindent We shall now state one more consequence of applying resolution of singularities to the Denjoy-Carleman quasianalytic class $C^\R_M ,$ the uniformization theorem, which is Corollary 5.14 on p.18 of \cite{bm} for the support of an ideal sheaf of finite type and also follows from the arguments in \cite{bmv} for the support of a quasicoherent ideal sheaf:

\smallskip
\newtheorem{unif}[quasian]{Theorem (Uniformization Theorem)}
\begin{unif}
Let $X$ be a manifold of class $C^\R_M,$ and let $Y$ be the support of any quasicoherent ideal sheaf that is a subsheaf of ${\mathfrak G}_X^{C^\R_M},$ the sheaf of germs of functions of class $C^\R_M$ defined at the points of $X.$ There exist a manifold $N$ and a proper $C^\R_M$-mapping $\varphi: N \rightarrow X$ such that $\varphi(N)=Y.$
\label{uniformization}
\end{unif}

\smallskip\noindent By Remark 5.15 in \cite{bm}, Theorem~\ref{uniformization} holds with $\dim N= \dim Y.$ We shall employ this conclusion for germs of varieties of class $C_M$ shortly, but before we can do so, we must discuss some important results that follow from Corollaries~\ref{topnoetherianity} and ~\ref{lojasiewicz}. I am indebted to \cite{pieroni} for the observation contained in the first of these:

\smallskip
\newtheorem{allfinitetype}[quasian]{Proposition}
\begin{allfinitetype}
Let $\I$ be any ideal in the ring of germs $C_M (x_0)$ for $x_0 \in \R^n$ or $x_0 \in \C^n$ or in $ C_{M,X}(x_0)$  for $X$ a manifold of class $C_M$ and $x_0 \in X,$
then there exists a subideal of finite type $\J \subset \I$ such
that $\forall \, f \in \I,$ $f \equiv 0$ on $\V(\J).$
\label{finitetypeprop}
\end{allfinitetype}

\smallskip\noindent {\bf Remark:} The paragraph preceding Definition~\ref{quasicoherent} explains how we choose a representative for $\V(\I)$ if $\I$ is finitely generated. As a consequence of this proposition, we choose $\V(\J)$ as a representative of $\V(\I),$ if $\I$ is infinitely generated.

\smallskip\noindent {\bf Proof:} The completion of the ring of germs $C_M (x_0)$ is the ring of formal power series $\F_{x_0}$ at $x_0,$ which is a Noetherian ring. Therefore, while $\I$ may not be finitely generated, $\I \cdot \F_{x_0}$ is finitely generated as an ideal of $\F_{x_0}.$ Let $f_1, \dots, f_k \in \I$ be generators for $\I \cdot \F_{x_0}.$ Due to quasianalyticity, $f$ is identified with its Taylor series at $x_0,$ so we can view it as an element of $\I \cdot  \F_{x_0},$ but then $f = \sum_{i=0}^k G_i f_i,$ where $G_i \in \F_{x_0}.$ Therefore, the Taylor series of $f$ at $x_0$ vanishes formally near $x_0$ on every $C_M$ curve in $\V( (f_1, \dots, f_k))$ with initial point $x_0.$ Quasianalyticity then implies that $f \equiv 0$ on $\V( (f_1, \dots, f_k)).$ We set $\J= (f_1, \dots, f_k).$ The reader should note that the argument given for the statement $\forall \, f \in \I,$ $f \equiv 0$ on $\V(\J)$ also shows that choosing another set of generators for $\I \cdot \F_{x_0}$ does not modify $\V(\J)$ apart from potentially shrinking the common neighborhood around $x_0$ on which all generators of $\I \cdot \F_{x_0}$ are defined. \qed

\smallskip\noindent The previous result and the {\L}ojasiewicz
inequalities, Corollary~\ref{lojasiewicz} imply a {\L}ojasiewicz
type Nullstellensatz. Before we can state this Nullstellensatz, we
have to specify which notion of radical we will employ:

\smallskip
\newtheorem{realraddef}[quasian]{Definition}
\begin{realraddef}
Let $C_M (x_0)$ for $x_0 \in \R^n$ or $x_0 \in \C^n$ be any ring of Denjoy-Carleman quasianalytic germs, and let $J \subset C_M (x_0)$ be an ideal, then the real radical of $J$
denoted by $\sqrt[\R]{J}$ is the ideal $$\left\{ g \in C_M(x_0) \: \big| \: \exists \, U \text{open}, \: U \ni x_0, \: \exists \, m \in \Q^+ , \: \exists \, f \in J, \: \: \text{such that} \: |g|^{m} \leq |f | \: \text{on} \: U \right\}.$$ Similarly, if $J \subset C_{M,X}(x_0)$ for $X$ a manifold of class $C_M$ and $x_0 \in X,$ $$\sqrt[\R]{J}=\left\{ g \in C_{M,X}(x_0) \: \big| \: \exists \, U \text{open}, \: U \ni x_0, \: \exists \, m \in \Q^+ , \: \exists \, f \in J, \: \: \text{such that} \: |g|^{m} \leq |f | \: \text{on} \: U \cap X \right\}.$$ 
\label{realrad}
\end{realraddef}

\smallskip
\newtheorem{nss}[quasian]{Theorem ({\L}ojasiewicz Nullstellensatz)}
\begin{nss}
Let $\I=(f_1, \dots, f_p)$ be any ideal
of finite type in $C_M (x_0)$ or in $C_{M,X}(x_0)$ for $X$ a manifold of class $C_M$ and $x_0 \in X.$ Let $\V(\I)$ be the germ of a variety corresponding to $\I,$ and let $\I(\V(\I))$ be the
ideal of germs of functions in $C_M (x_0)$ or in $C_{M,X}(x_0)$ vanishing on $\V(\I).$
Then $$\sqrt[\R]{\I}=\I(\V(\I)).$$ \label{nullstellensatz}
\end{nss}

\smallskip\noindent {\bf Proof:} The inclusion $\sqrt[\R]{\I}\subset\I(\V(\I))$
is clear. We only have to prove the reverse inclusion. Without loss of generality, we only prove the result for the ring $C_M (x_0)$ noting that the case of $C_{M,X} (x_0)$ is similar, only requiring that all open sets are intersected with $X.$ Consider
any $h \in \I(\V(\I)).$ Since $h \in C_M(x_0),$ there exists an open set $U_h \ni x_0$ such that  $h \in C_M(U_h).$ Similarly, for each $i,$ $1 \leq i \leq p,$ there exists $U_i \ni x_0$ open such that $f_i \in C_M(U_i).$ Set $\tilde U = U_h \, \cap \, U_1 \, \cap \dots \cap \, U_p$ and $g=|f_1|^2+\cdots+|f_p|^2.$ Clearly, $\tilde U$ is open, and $g,h \in C_M(\tilde U).$  Up to potentially shrinking $\tilde U,$ $$\{ x \in \tilde U \: \big| \: g(x)=0 \} \subseteq \{ x
\in \tilde U \: \big| \: h(x)=0 \}$$ and $\{x_0\}$ is a compact set, so we apply part (I) of
Corollary~\ref{lojasiewicz} to conclude
that there exist $C, \alpha >0$ and an open neighborhood $U_{x_0} \subset \tilde U$ with $x_0 \in U_{x_0}$ such that
$$|g(x)| \geq C \, |h(x)|^{\, \alpha}$$ for all $x \in U_{x_0}.$ The exponent
$\alpha$ is rational. This precisely means that $h \in
\sqrt[\R]{\I}$ as needed. \qed

\smallskip\noindent We now use the {\L}ojasiewicz
inequalities to deduce that the germ of a variety corresponding to
an ideal $\I$ in a Denjoy-Carleman quasianalytic ring of germs $C_M (x_0)$ or $C_{M,X}(x_0)$ for $X$ a manifold of class $C_M$ and $x_0 \in X$ must have an open and dense set of smooth points. First,
let us define a smooth point of an affine variety:

\smallskip
\newtheorem{smoothpt}[quasian]{Definition}
\begin{smoothpt}
Let $\V= \V (\I)$ be a variety corresponding to an ideal $\I$
in $C_M (U)$ (or in $C_{M,X}(U)$ for $X$ a manifold of class $C_M$), where $U$ is an open set in $\C^n$ (or in $X$). $x' \in \V$ is a smooth point of $\V$ if there
exist a neighborhood $\tilde U \ni x'$ (open in the topology of $X$ in the case of a variety corresponding to an ideal in $C_{M,X}(U)$) and smooth functions $f_1, \dots, f_s$ such that $$\V \cap \tilde U = \{x \in \tilde U \: | \: f_1(x)= \cdots = f_s(x)=0 \} $$ is a $\smooth$
submanifold of $\C^n$ (or of $X$) of codimension $s,$ i.e. $d f_1
\wedge \cdots \wedge
d  f_s (x) \neq 0$ for all $x \in \tilde U.$ \label{smoothptdef}
\end{smoothpt}

\smallskip\noindent {\bf Remark:} The definition implies that the set of smooth points
of a variety $\V$ is open.

\bigskip\noindent For the ring of holomorphic functions $\oka$ and
for the ring of $\C$-valued real analytic functions $\ra,$ it is
known that any affine variety has an open and dense set of smooth
points. For the ring of smooth functions $\smooth,$ this same
statement is not always true but was proven by Ren{\'e} Thom in
\cite{thom} for $\V = \V(\I)$ under the hypothesis that the ideal
$\I \subset \smooth$ is {\L}ojasiewicz. An ideal $\I \subset
\smooth$ is called {\L}ojasiewicz if it is of finite type and its
generators satisfy the {\L}ojasiewicz inequality with respect to
distance, which is the content of part (II) of
Corollary~\ref{lojasiewicz} for all elements of the ring $C_M
(U).$ Given any ideal $\I $ in $C_M (x_0)$ or in $C_{M,X}(x_0),$
Proposition~\ref{finitetypeprop} guarantees that $\V(\I)$ can be
presented as the variety $\V(\J)$ corresponding to a subideal $\J$ of $\I$ of finite type.
This result combined with part (II) of Corollary~\ref{lojasiewicz}
leads us to expect that indeed any germ of a variety corresponding to
an ideal in $C_M (x_0)$ or in $C_{M,X}(x_0)$ must have an open and dense set of smooth points:

\smallskip
\newtheorem{smoothptsdense}[quasian]{Proposition}
\begin{smoothptsdense}
Let $\I$ be any ideal in $C_M (x_0)$ for $x_0 \in \R^n$ or $x_0 \in \C^n$ or in $C_{M,X}(x_0)$ for $X$ a manifold of class $C_M$ and $x_0 \in X.$ Let
$\V(\I)$ be any representative of the variety corresponding to
$\I,$ then $\V(\I)$ has an open and dense set of smooth points.
\label{smoothptsprop}
\end{smoothptsdense}

\smallskip\noindent {\bf Proof:} We know the set of smooth points
of $\V(\I)$ is open, so we only have to show that it is dense. Just like in the proof of Theorem~\ref{nullstellensatz}, without loss of generality we will prove this result only for $C_M (x_0)$ noting once again that for $C_{M,X}(x_0)$ all open sets are intersected with $X.$ We apply
Proposition~\ref{finitetypeprop} to conclude that there exists a
subideal of finite type $\J= (f_1, \dots, f_p) \subset \I$ such
that $\V(\J)$ can be taken as a representative of $\V(\I).$ Set $f=|f_1|^2+\cdots+|f_p|^2.$ For each $i,$ $1 \leq i \leq p,$ there exists $U_i \ni x_0$ open such that $f_i \in C_M(U_i).$ Set $\tilde U = U_1 \, \cap \dots \cap \, U_p.$ Clearly, $\tilde U$ is open, $f \in C_M(\tilde U),$ and $\V(\I) \, \cap \, \tilde U = \{ x \in \tilde U \: \big| \: f(x)=0 \}.$
By part (II) of Corollary~\ref{lojasiewicz}, the {\L}ojasiewicz
inequality with respect to distance holds for $f$ on $\tilde U$ up to potentially shrinking $\tilde U,$ so the
proof of Rene Thom in \cite{thom} or the even simpler proof of the
same result given by Jean-Claude Tougeron in \cite{tougeron}
(proof of Proposition $4.6$ in subsection $V.4$) apply verbatim. The reader should note that given any subideal $\J$ of $\I$ satisfying the conclusion of Proposition~\ref{finitetypeprop}, this proof applies to show that $\V(\J)$ taken as a representative for $\V(\I)$ has an open and dense set of smooth points.
\qed

\smallskip\noindent {\bf Remark:} Thom's theorem does not apply
directly to this case because we are considering here an ideal
$\I$ of elements in the ring $C_M (x_0)$ or $C_{M,X}(x_0)$ and not a \L ojasiewicz ideal in the ring of smooth functions.

\medskip\noindent In fact, the Uniformization Theorem, Theorem~\ref{uniformization}, implies an even stronger statement than Proposition~\ref{smoothptsprop} for the point $x_0$ itself:

\smallskip
\newtheorem{regularptsdense}[quasian]{Proposition}
\begin{regularptsdense}
Let $\I$ be any ideal in $C_M (x_0)$ for $x_0 \in \R^n$ or $x_0 \in \C^n$ or in $C_{M,X}(x_0)$ for $X$ a manifold of class $C_M$ and $x_0 \in X.$ Let
$\V(\I)$ be any representative of the variety corresponding to
$\I,$ let $\I=\I (\V(\I)),$ and let $\V(\I)$ not consist of an isolated point ($\V(\I) \neq \{x_0\}),$ then $x_0$ is in the closure of the set of regular points of $\V(\I)$ of maximal dimension. \label{regularptsprop}
\end{regularptsdense}

\smallskip\noindent {\bf Proof:} Once again, we will prove this result only for $C_M (x_0)$ for $x_0 \in \R^n.$ Just like in the proof of Proposition~\ref{smoothptsprop}, we first apply Proposition~\ref{finitetypeprop} to obtain a subideal of finite type $\J= (f_1, \dots, f_p) \subset \I$ such that $\V(\J)$ can be taken as a representative for $\V(\I).$ Without loss of generality, $f_i \in C_M^\R(x_0)$ for all $1\leq i \leq p.$ Furthermore, since $\I=\I( \V(\I)),$ we can choose $f_1, \dots, f_p$ in such a way that $\V(\I)$ gets the reduced structure. Consider $f_1.$ Without loss of generality, we can choose $f_1$ such that there exists $j,$ $1 \leq j \leq n$ and $\{x^{(k)}\} \subset \V(\I)$ a sequence of points in $\V(\I)$ converging to $x_0$ such that $\frac{\partial f_1}{\partial x_j} (x^{(k)}) \neq 0.$ This is true because if it were not the case, then each derivative of $f_1$ $\frac{\partial f_1}{\partial x_j}$ for $1 \leq j \leq n$ would be identically zero on $\V(\I)$ in a neighborhood $U_j$ of $x_0.$ The ring $C_M^\R (x_0)$ is closed under differentiation, so $\frac{\partial f_1}{\partial x_j} \in C_M^\R(U_j),$ and also  $\I=\I (\V(\I)),$ so $\frac{\partial f_1}{\partial x_j} \in \I.$ We can thus keep looking at higher order derivatives of $f_1$ shrinking the neighborhood around $x_0$ as necessary until we find one with the property stated above. Note that this process finishes in finitely many steps because every function in $C_M^\R(x_0)$ has a finite order of vanishing at $x_0.$ 
We carry out the same procedure for $f_2$ and so on until we have chosen a subideal of finite type $\J= (f_1, \dots, f_p) \subset \I$ that gives the reduced structure on $\V(\I).$ Now, for each $i,$ $1 \leq i \leq p,$ there exists $U_i \ni x_0$ open such that $f_i \in C^\R_M(U_i).$ Set $\tilde U = U_1 \, \cap \dots \cap \, U_p.$ Note that  $f_i \in C^\R_M(\tilde U)$ for each $i.$ Consider $\V(\J)$ on $\tilde U.$ It is the support of an ideal sheaf of finite type with the reduced structure and no isolated points, so the Uniformization Theorem, Theorem~\ref{uniformization}, implies that there exist a manifold $N$ and a proper desingularization mapping $\varphi: N \rightarrow \R^n$ of class $C^\R_M$ such that $\varphi(N)= \V(\I) \cap \tilde U,$ so $x_0$ is in the closure of regular points of maximal dimension of $\V(\I)$ as needed.  \qed

\medskip\noindent {\bf Remarks:} 

\smallskip \noindent (1) In the proof of Theorem~\ref{maintheorem}, we are not concerned whether the sequence of regular points converging to $x_0$ consists of regular points of maximal dimension. Just being regular points of the variety will suffice.

\smallskip\noindent (2) The way we chose $f_1$ in this proof is essentially the way Bierstone and Milman bypass in \cite{bm} the absence of the Weierstrass Preparation Theorem for the Denjoy-Carleman classes of functions. The interested reader should consult Remark 7.10 of \cite{bm} regarding this matter, reduced structures, and alternative ways of running the desingularization algorithm.

\smallskip\noindent (3) For any subideal $\J$ of $\I$ satisfying the conclusion of Proposition~\ref{finitetypeprop}, the same argument applies to obtain the same conclusion for any representative $\V(\J)$ of $\V(\I).$

\bigskip\noindent Finally, we can prove Theorem~\ref{radaccprop}:

\smallskip\noindent {\bf Proof of Theorem~\ref{radaccprop}:} Let $\I_1 \subset \I_2 \subset \I_3 \subset \dots$ be an increasing chain of ideals. Consider the increasing chain of ideals $\sqrt[\R]{\I_1} \subset \sqrt[\R]{\I_2} \subset \sqrt[\R]{\I_3} \subset \dots$ By Proposition~\ref{finitetypeprop}, for each $j \geq 1,$ we know that there exists a subideal of finite type $\J_j \subset  \sqrt[\R]{\I_j}$ such that $\V(\J_j)$ can be taken as a representative for $\V(\sqrt[\R]{\I_j}).$ We apply next the {\L}ojasiewicz Nullstellensatz, Proposition~\ref{nullstellensatz}, to conclude that $\I(\V(\sqrt[\R]{\I_j}))=\I(\V(\J_j))=\sqrt[\R]{\J_j}$ for all $j \geq 1$ and both $\I(\V(\sqrt[\R]{\I_j}))$ and $\I(\V(\J_j))$ are ideals in $C_M (x_0).$ Note that we needed a finitely generated ideal in order to be able to apply the Nullstellensatz. By Topological Noetherianity, Corollary \ref{topnoetherianity}, the decreasing chain of germs of varieties $$\V(\sqrt[\R]{\I_1}) \supset \V(\sqrt[\R]{\I_2}) \supset \dots$$ stabilizes, namely there exists some $k \in \N$ such that for all $j \geq k,$ $\V(\sqrt[\R]{\I_j})=\V(\sqrt[\R]{\I_k}).$ The largest ideal of functions, however, that vanishes on $\V(\sqrt[\R]{\I_j})$ is $\I(\V(\sqrt[\R]{\I_j})).$ Therefore, the increasing chain of ideals $$\I(\V(\sqrt[\R]{\I_1})) \subset \I(\V(\sqrt[\R]{\I_2})) \subset \dots$$ stabilizes. We know $\sqrt[\R]{\J_j} \subset \sqrt[\R]{\I_j} \subset \I(\V(\sqrt[\R]{\I_j})),$ so from the application of the Nullstellensatz above, we conclude $\I(\V(\sqrt[\R]{\I_j}))=   \sqrt[\R]{\I_j}$ for every $j \geq 1,$ so $\sqrt[\R]{\I_1} \subset \sqrt[\R]{\I_2} \subset \sqrt[\R]{\I_3} \subset \dots$ stabilizes. \qed

\smallskip\noindent {\bf Remark:} It should be noted here that Theorem~\ref{radaccprop} also holds for an ascending chain of ideals $\I_1 \subset \I_2 \subset \cdots$ in $C_{M,X} (x_0),$ where $X$ is a manifold of class $C_M$ and $x_0 \in X,$ i.e. the ascending chain of radical ideals $\sqrt[\R]{\I_1} \subset \sqrt[\R]{\I_2} \subset \cdots$ stabilizes. We will use this observation in the proof of Theorem~\ref{maintheorem}.

\bigskip\bigskip

\section{Subellipticity of the $\dbar$-Neumann Problem and the
Kohn Algorithm} \label{subkohnalg}

\bigskip

We shall give here only a brief outline of the properties of
subelliptic multipliers and the Kohn algorithm. The interested
reader should consult Kohn's paper \cite{kohnacta} for
the set-up of the $\dbar$-Neumann problem as well as the proofs of
the results cited in this section. We first give Kohn's
characterization of subellipticity of the $\dbar$-Neumann problem
on $(p,q)$ forms:

\smallskip
\newtheorem{subellgain}{Definition}[section]
\begin{subellgain}
Let $\Omega$ be a domain in $\C^n$ and let $x_0 \in
\overline{\Omega}.$ The $\dbar$-Neumann problem on $\Omega$ for
$(p,q)$ forms is said to be subelliptic at $x_0$ if there exist a
neighborhood $U$ of $x_0$ and constants $C, \epsilon > 0$ such
that
\begin{equation}
||\varphi\, ||_{\, \epsilon}^2 \leq C \, ( \, ||\,\dbar \,
\varphi\,||^2_{\, 0} + ||\, \ad  \varphi \,||^2_{\, 0} +
||\,\varphi \,||^2_{\, 0} \,) \label{subellestproblem}
\end{equation}
for all $(p,q)$ forms $\varphi \in \smoothc (U) \cap Dom (\ad),$
where $||\, \cdot \, ||_{\, \epsilon}$ is the Sobolev norm of
order $\epsilon$ and $||\, \cdot \,||_{\, 0}$ is the $L^2$ norm.
\label{subellgaindef}
\end{subellgain}

\smallskip\noindent The boundary condition $\varphi \in Dom(\ad)$ is responsible for the non-ellipticity of the $\dbar$-Neumann problem. If the point $x_0$ is inside the domain
$\Omega,$ then automatically estimate \eqref{subellestproblem} holds
at $x_0$ with the largest possible value $\epsilon = 1.$ The problem is
thus elliptic rather than subelliptic inside; see \cite{dbneumann1} and \cite{dbneumann2} for strongly
pseudoconvex domains and \cite{hormest}
and \cite{kohnweakpsc} for
pseudoconvex domains. Therefore, subellipticity only needs to be
studied on the boundary of the domain $b \Omega.$

We shall now give Kohn's definition of a subelliptic multiplier:

\smallskip
\newtheorem{subellmult}[subellgain]{Definition}
\begin{subellmult}
Let $\Omega$ be a domain in $\C^n$ and let $x_0 \in
\overline{\Omega}.$ A $\smooth$ function $f$ is called a
subelliptic multiplier at $x_0$ for the $\bar\partial$-Neumann
problem on $\Omega$ if there exist a neighborhood $U$ of $x_0$ and
constants $C, \epsilon > 0$ such that
\begin{equation}
||\, f \varphi\, ||_{\, \epsilon}^2 \leq C \, ( \, ||\,\dbar \,
\varphi\,||^2_{\, 0} + ||\, \ad  \varphi \,||^2_{\, 0} +
||\,\varphi \,||^2_{\, 0} \,) \label{subellest}
\end{equation}
for all $(p,q)$ forms $\varphi \in \smoothc (U) \cap Dom (\ad).$
We will denote by $I^q (x_0)$ the set of all subelliptic
multipliers at $x_0.$ \label{subellmultdef}
\end{subellmult}

\smallskip\noindent The notation $I^q (x_0)$ for subelliptic
multipliers at $x_0$ for $(p,q)$ forms drops reference to $p,$ the
holomorphic part of such forms, which is irrelevant in the
$\dbar$-Neumann problem.

\medskip\noindent {\bf Remarks:}

\smallskip \noindent (1) Clearly, if there is a subelliptic multiplier $f \in I^q (x_0)$ such that $f(x_0) \neq 0$, then a
subelliptic estimate holds at $x_0$ for the $\dbar$-Neumann
problem.

\smallskip\noindent (2) Estimate \eqref{subellest} holds for the largest possible value $\epsilon=1$ if $f = 0$ on $U \cap b \Omega.$ An example of such $f$ is $r,$ the defining function of the domain $\Omega.$

\smallskip\noindent (3) If $x_0 \in b \Omega,$ the highest
possible gain in regularity in the $\dbar$-Neumann problem given
by $\epsilon$ in estimate~\eqref{subellestproblem} is
$\frac{1}{2}$ for $\Omega$ a strongly pseudoconvex domain; see
\cite{dbneumann1} and \cite{dbneumann2}.

\smallskip\noindent (4) The reader should note that Definition~\ref{subellmultdef} yields a sheaf of germs of subelliptic multipliers, which is precisely what is needed to capture the subellipticity of the $\dbar$-Neumann problem. The underlying ring is $\smooth,$ so pointwise $I^q (x_0) \subset \smooth_{b \Omega} (x_0),$ the ring of germs of smooth functions on $b \Omega$ at $x_0 \in b \Omega.$ Globally, we will denote by $\E_{b \Omega}$ the sheaf of germs of smooth functions defined on the boundary of the domain $b \Omega$ as to match the notation for the sheaf of germs of smooth functions used by Malgrange and his school. Hence  $\smooth_{b \Omega} (x_0)=\E_{b \Omega}(x_0).$

\medskip\noindent Kohn's Theorem $1.21$ in \cite{kohnacta} encapsulates the properties of
subelliptic multipliers he proved in his paper. These motivate the
way he sets up his algorithm, which determines whether or not the
$\dbar$-Neumann problem is subelliptic:

\medskip
\newtheorem{subellcor}[subellgain]{Theorem}
\begin{subellcor}
If $\Omega$ is pseudoconvex with a $\smooth$ boundary and if $x_0
\in \overline{\Omega},$ then we have:
\begin{enumerate}
\item[(a)] $I^q (x_0)$ is an ideal.
\item[(b)] $I^q (x_0) = \sqrt[\R]{I^q (x_0)}.$
\item[(c)] If $r=0$ on $b \Omega,$ then $r \in I^q (x_0)$ and the
coefficients of $\partial r \wedge \dbar r \wedge (\partial \dbar
r)^{n-q}$ are in $I^q (x_0).$
\item[(d)] If $f_1, \dots, f_{n-q} \in I^q (x_0),$ then the
coefficients of $\partial f_1 \wedge \dots \wedge \partial f_j
\wedge \partial r \wedge \dbar r \wedge (\partial \dbar
r)^{n-q-j}$ are in $I^q (x_0),$ for $j \leq
n-q.$\label{subellpropcor}
\end{enumerate}
Here $\sqrt[\R]{I^q (x_0)}$ is computed in the ring of germs $\smooth_{b \Omega} (x_0)=\E_{b \Omega}(x_0),$ i.e. $\sqrt[\R]{I^q (x_0)}$ consists of all $g \in \smooth_{b \Omega} (x_0)$ such that there exist an open set $U_{x_0}$ containing $x_0,$ some $\alpha \in \Q^+ ,$ and $f \in I^q(x_0)$ satisfying $|g(x)| \leq |f(x)|^\alpha$ for all $x \in U_{x_0} \cap b \Omega.$
\end{subellcor}

\medskip\noindent {\bf The Kohn Algorithm:}

\medskip\noindent {\bf Step 1}  $$I^q_1(x_0) = \sqrt[\R]{(\, r,\,
\text{coeff}\{\partial r \wedge \dbar r \wedge (\partial \dbar
r)^{n-q}\}\, )}$$

\medskip\noindent {\bf Step (k+1)} $$I^q_{k+1} (x_0) = \sqrt[\R]{(\, I^q_k
(x_0),\, A^q_k (x_0)\, )},$$ where $$A^q_k (x_0)= \text{coeff}\{\partial
f_1 \wedge \dots \wedge \partial f_j \wedge \partial r \wedge
\dbar r \wedge (\partial \dbar r)^{n-q-j}\}$$ for $f_1, \dots, f_j
\in I^q_k (x_0)$ and $j \leq n-q.$ Notation $( \, \cdot \, )$
denotes the ideal generated by the functions inside the
parentheses in the ring $\smooth_{b \Omega} (x_0),$ and $\text{coeff}\{\partial r \wedge \dbar r \wedge (\partial\dbar r)^{n-q}\}$ is the determinant of a minor of the determinant of the Levi form for $q>1$ and equals the determinant of the Levi form itself for $q=1.$ It is easily seen that $I^q_k (x_0) \subset I^q (x_0)$ at each step $k$ and
that the Kohn algorithm generates an increasing chain of ideals
$$I^q_1 (x_0) \subset I^q_2 (x_0) \subset \cdots$$ in the ring of germs $\smooth_{b \Omega} (x_0)=\E_{b \Omega}(x_0).$ We will denote by $I^q_k$ the subelliptic multipliers at step $k$ of the Kohn algorithm on $b \Omega,$ which is a sheaf of ideals that is a subsheaf of $=\E_{b \Omega}.$ We will employ the
following notation pertaining to the varieties corresponding to
ideals of multipliers just as Kohn does in \cite{kohnacta}:
$$\V_k^q (x_0) = \V (I^q_k (x_0))$$ defined at each $x_0 \in b \Omega$ as well as the global notation $$\V_k^q = \V (I^q_k),$$ where $$\V (I^q_k)= supp \: \frac{ \E_{b \Omega}}{I^q_k},$$ which is consistent with Section~\ref{DCdescription}.

To understand the behavior of the sheaves $I^q_k,$ we need to introduce the property of being quasi-flasque, a notion defined by Tougeron in \cite{tougeronqf}:

\medskip
\newtheorem{qfdef}[subellgain]{Definition}
\begin{qfdef}
Consider an open set $\tilde U$ in $\R^m$ and
$\E$ the sheaf of $\smooth$ germs on $\tilde U.$ A sheaf $\M$
of $\E$-modules is quasi-flasque if for every open set $U
\subset \tilde U$ the canonical homomorphism
$$\M (\tilde U) \otimes_{\E(\tilde U)} \E (U) \longrightarrow \M
(U)$$ is an isomorphism.
\end{qfdef}

As proven in \cite{andreeaqf}, the sheaves of multipliers in the Kohn algorithm $I^q_k$ are quasi-flasque:

\medskip
\newtheorem{qfsteps}[subellgain]{Theorem}
\begin{qfsteps}
The sheaf $I^q_k$ satisfies the following two properties:
\begin{enumerate}
\item[(a)] $I^q_k$ is quasi-flasque;
\item[(b)] Let $x_0$ be any point in $b\Omega.$ If sections $s_j \in I^q_k(b \Omega)$ generate $I^q_k(x_0)$ for $j \in J,$ $J$ an indexing set, then $s_j$ also generate $I^q_k(x)$ for $x$ sufficiently close to $x_0.$
\end{enumerate}
\label{quasiflasquesteps}
\end{qfsteps}

In other words, the sheaves $I^q_k$ are quasicoherent. Note that due to the form of estimate \eqref{subellmultdef}, a subelliptic multiplier on some neighborhood $U$ can be extended to a global multiplier on $b \Omega$ by trivially extending it to be identically zero outside of $U;$ see \cite{andreeaqf}. By the way, the fact that any quasi-flasque sheaf has the property in Theorem~\ref{quasiflasquesteps} (b) is Tougeron's Proposition 6.4 from section V.6 of \cite{tougeron}.

We will now recall from \cite{kohnacta} the definition of the
Zariski tangent space to an ideal and to a variety, which are
crucial in testing the progress of the Kohn algorithm. We will
state everything in terms of ideals in the ring $\smooth_{b \Omega} (x_0),$ which
is the underlying pointwise ring in our construction, but we note that these definitions are identical for ideals in $\ra_{b \Omega}(x_0)$ or $C_{M, b \Omega}(x_0).$

\smallskip
\newtheorem{zariskidef}[subellgain]{Definition}
\begin{zariskidef}
Let $\I$ be an ideal in $\smooth_{b \Omega} (x_0),$ and let $\V (\I)$ be a germ of the variety corresponding to $\I.$ We define $Z^{\,
1,0}_x (\I)$ the Zariski tangent space of $\I$ at some $x \in \V(\I)$ to be
$$Z^{\, 1,0}_x (\I) = \{ \, L \in T^{\, 1,0}_x (b \Omega) \: | \: L(f)=0 \:\:
\forall \: f \in \I \, \},$$ where $T^{\, 1,0}_x (b \Omega)$ is the
$(1,0)$ tangent space to $b \Omega \subset \C^n$ at $x.$ If $\V$ is a germ of a
variety at $x_0$ and $x \in \V,$ then we define
$$Z^{\, 1,0}_x (\V) = Z^{\, 1,0}_x (\I_x (\V)), $$ where $\I_x (\V)$
is the ideal of all germs of functions in $\smooth_{b \Omega} (x)$ vanishing on $\V.$
\end{zariskidef}

\medskip\noindent The next lemma is Lemma $6.10$ of \cite{kohnacta}
that relates $Z^{\, 1,0}_x (\I)$ with $Z^{\, 1,0}_x (\V(\I)).$ We state it here for ideals in $\smooth_{b \Omega} (x_0),$ but it is also true for ideals in $\ra_{b \Omega}(x_0)$ or $C_{M, b \Omega}(x_0):$

\smallskip
\newtheorem{zariskilemma}[subellgain]{Lemma}
\begin{zariskilemma}
\label{zariski} Let $\I$ be an ideal in $\smooth_{b \Omega} (x_0),$ then
\begin{equation}
Z^{\, 1,0}_x (\V(\I)) \subset Z^{\,1,0}_x (\I). \label{zarexpr}
\end{equation}
Equality holds in ~\eqref{zarexpr} if the ideal $\I_x (\V(\I))$ is generated by elements of $\I.$ 
\end{zariskilemma}

\smallskip\noindent Equality also holds in ~\eqref{zarexpr} if $x$ is a regular point of the variety $\V(\I).$ To prove this fact, we recall here an elementary lemma, which can be found on p.21 of \cite{boggess}:

\smallskip
\newtheorem{manifoldnsslemma}[subellgain]{Lemma}
\begin{manifoldnsslemma}
Suppose $Y$ is an $l$-dimensional smooth submanifold of $\R^N$ given by $Y=\{\rho_1=\cdots=\rho_{N-l}=0\}$ with $d \rho_1 \wedge \dots \wedge d\rho_{N-l} \neq 0$ on $Y.$ Suppose $f: \R^N \rightarrow \R$ is a smooth function that vanishes on $Y.$ Then there are smooth functions $\alpha_1, \dots, \alpha_{N-l}$ defined near $Y$ so that $$f=\sum_{j=1}^{N-l} \: \alpha_j \, \rho_j \quad \text{near} \quad Y.$$
\label{manifoldnss}
\end{manifoldnsslemma}

\smallskip\noindent Now we can prove our claim:

\smallskip
\newtheorem{zariskireglemma}[subellgain]{Lemma}
\begin{zariskireglemma}
Let $\I$ be an ideal in $\smooth_{b \Omega} (x_0),$ $\ra_{b \Omega}(x_0),$ or $C_{M, b \Omega}(x_0).$ For any representative of $\V(\I),$ if $x$ is a regular point of $\V(\I),$ then $$Z^{\, 1,0}_x (\V(\I)) =Z^{\,1,0}_x (\I).$$
\label{zariskireg}
\end{zariskireglemma}

\smallskip\noindent {\bf Proof:} By Lemma~\ref{zariski}, $Z^{\, 1,0}_x (\V(\I)) \subset Z^{\,1,0}_x (\I),$ so we only need to prove the reverse inclusion. Since $x$ is a regular point of $\V(\I),$ there exist a neighborhood $U \ni x$ and functions $\rho_1, \dots, \rho_{N-l} \in \I$ such that $\V(\I) \cap U = \{\rho_1=\cdots=\rho_{N-l}=0\}$ with $d \rho_1 \wedge \dots \wedge d\rho_{N-l} \neq 0$ on all of $U.$ Consider an arbitrary $L \in Z^{\,1,0}_x (\I).$ For any $f \in \I_x (\V(\I)),$ by Lemma~\ref{manifoldnss}, there exist smooth functions $\alpha_1, \dots, \alpha_{N-l}$ such that on $U$ or a shrinking of $U,$  $f=\sum_{j=1}^{N-l} \: \alpha_j \, \rho_j.$ $$L f = \sum_{j=1}^{N-l} \: L(\alpha_j \, \rho_j) =  \sum_{j=1}^{N-l} \: \left(L(\alpha_j) \, \rho_j+ \alpha_j \, L(\rho_j)\right)$$ Since $L \in Z^{\,1,0}_x (\I)$ and $\rho_j \in \I$ for each $j,$ $L(\rho_j)(x)=0$ for each $j.$ Since $x \in \V(\I) \cap U,$ $\rho_j (x)=0$ for each $j.$ We conclude that $Lf (x)=0,$ which implies $L \in Z^{\, 1,0}_x (\V(\I))$ as needed. \qed

\bigskip\noindent We define $$\nx = \{ \, L \in T^{\, 1,0}_x (b \Omega) \: |
\: \langle \, (\partial \dbar r)_x \, , \, L \wedge \bar L \,
\rangle = 0 \, \},$$ which is the subspace of $T^{\, 1,0}_x (b
\Omega)$ consisting of the directions in which the Levi form
vanishes. With this notation, we can define the holomorphic
dimension of a variety lying in the boundary of the domain
$\Omega.$

\smallskip
\newtheorem{holdimdef}[subellgain]{Definition}
\begin{holdimdef}
Let $\V$ be a variety such that $\V \subset b \Omega.$ Let $x \in \V.$ We define the holomorphic dimension of $\V$ at $x$ by $$hol. \: dim_x \, (\V) =  \dim Z^{\, 1,0}_x (\V) \cap \nx.$$ We define the holomorphic dimension of the entire variety $\V$ by $$hol. \: dim \, (\V) = \min_{x
\in \V} \dim Z^{\, 1,0}_x (\V) \cap \nx.$$
\end{holdimdef}

\smallskip\noindent We conclude this section with  Proposition $6.12$ of \cite{kohnacta} that explores the
implications of the Kohn algorithm not advancing at a particular
point:

\smallskip
\newtheorem{stillthere}[subellgain]{Proposition}
\begin{stillthere}
If $x \in \V_k^q(x_0),$ then $$x \in \V_{k+1}^q (x_0)  \quad
\Leftrightarrow \quad \dim (Z^{\, 1,0}_x (I^q_k (x_0)) \cap \nx)
\geq q.$$ \label{stilltherethm}
\end{stillthere}

\smallskip\noindent {\bf Remark:} The reader should note that Kohn's proof of this result carries through also for modifications of his algorithm such that only real analytic functions are considered (section 6 of \cite{kohnacta} and \cite{racase}) or only Denjoy-Carleman functions are considered (next section).

\bigskip\bigskip

\section{Equivalence of Types}
\label{equivtypepf}

\bigskip

We start by describing a modification of the Kohn algorithm on a domain whose defining function $r$ is of class $C_M$ that only produces germs of Denjoy-Carleman functions in the same class. We assume the class $C_M$ is closed under differentiation. As seen in Section~\ref{subkohnalg}, the Kohn algorithm involves only differentiations, taking of real
radicals, addition, and multiplication starting with $r.$ Addition and multiplication appear not only in the generation of various ideals but also in the computation of all Jacobians involved.  In the modified Kohn algorithm, ideals are generated in the ring $C_{M, b \Omega} (x_0)$ and $\sqrt[\R]{\, \cdot \, }$ is computed in the sense of Definition~\ref{realrad} solely with elements in the ring $C_{M, b \Omega} (x_0)$ and not all of $\smooth_{b \Omega} (x_0).$ It is clear that the termination of this modified algorithm implies the termination of the original Kohn algorithm on germs. We shall denote by $\tilde I^q_k (x_0)$ the ideals of subelliptic multipliers given by this modified Kohn algorithm.

\medskip\noindent {\bf The Modified Kohn Algorithm:}

\medskip\noindent {\bf Step 1}  $$\tilde I^q_1(x_0) = \sqrt[\R]{(\, r,\,
\text{coeff}\{\partial r \wedge \dbar r \wedge (\partial \dbar
r)^{n-q}\}\, )_{C_{M, b \Omega} (x_0)}}$$

\medskip\noindent {\bf Step (k+1)} $$\tilde I^q_{k+1} (x_0) = \sqrt[\R]{(\, \tilde I^q_k
(x_0),\, \tilde A^q_k (x_0)\, )_{C_{M, b \Omega} (x_0)}},$$ where $$\tilde A^q_k (x_0)= \text{coeff}\{\partial
f_1 \wedge \dots \wedge \partial f_j \wedge \partial r \wedge
\dbar r \wedge (\partial \dbar r)^{n-q-j}\}$$ for $f_1, \dots, f_j
\in \tilde I^q_k (x_0)$ and $j \leq n-q.$ Note that the subscript $C_{M, b \Omega} (x_0)$ emphasizes ideals are generated in this ring instead of $\smooth_{b \Omega} (x_0).$ $C_{M, b \Omega} (x_0)$ is a ring, hence closed under addition and
multiplication, and also closed under differentiation by assumption. By definition of the modified Kohn algorithm, all real radicals of ideals stay inside $C_{M, b \Omega} (x_0).$ Therefore, we conclude $\tilde I^q_k (x_0) \subset C_{M, b \Omega} (x_0)$ at every step $k$ of the Kohn algorithm.

Consider now the sheaves $\tilde I^q_k$ of Denjoy-Carleman subelliptic multipliers at all the steps of the modified Kohn algorithm. We would like to show they are quasicoherent, but that is difficult to do directly. Instead, we will show below that we can relate the sheaves $\tilde I^q_k$ in a natural way to the sheaves of the full scale Kohn algorithm $I^q_k,$ which are quasicoherent.  Define $\widehat{\tilde I^q_k} = \tilde I^q_k \otimes_{\E} \E$ by generating a new sheaf over the ring of smooth functions; see \cite{malgrange} and \cite{tougeronqf} for instances where the idea of tensoring with $\E$ is used to understand coherence properties of sheaves in the real analytic case. First, let us reprove a local version of Lemma 3.1 from \cite{abn}, where all rings of functions are defined on an open set $U \subset \R^n:$

\medskip
\newtheorem{raapprox}{Lemma}[section]
\begin{raapprox}
\label{ra approx}
Let $U \subset \R^n$ be an open set. Let $\ra(U)$ be the ring of real analytic functions on $U,$ and let $\An(U)$ be any ring of functions on $U$ such that $\ra(U) \subset \An(U) \subset \stb(U).$ If $\gta$ is any ideal in $\An (U)$ and if $\hat \gta$ is the ideal generated by the elements of $\gta$ in $\stb(U),$ then $  {\sqrt[\R]{\gta}}  =\sqrt[\R]{\hat\gta} \cap \An(U),$ where the real radical is computed on the entire open set $U.$
\end{raapprox}

\smallskip\noindent {\bf Remark:} The reader should note that $\gta$ is not assumed to be finitely generated.

\medskip\noindent {\bf Proof:} The inclusion  $  {\sqrt[\R]{\gta}}  \subset \sqrt[\R]{\hat\gta} \cap \An(U)$ is obvious. Thus, we only need to prove the reverse inclusion. Assume that $g \in \sqrt[\R]{\hat\gta} \cap \An(U).$ By the definition of the real radical, $\exists \,\tilde f \in \hat\gta$ and $\exists \, \alpha \in \Q^+$ such that $|g(x)|^\alpha \leq | \tilde f(x)|$ for every $x \in U.$ Note that the real radical here is computed on the whole open set $U$ unlike the definition for germs given in \ref{realrad}. Since $\alpha \in \Q^+,$ we can write $\displaystyle \alpha=\frac{p}{q}$ for $p, q \in \N^*.$ Therefore, $|g(x)|^p \leq | \tilde f(x)|^q$ for every $x \in U.$  We square both sides to obtain that $|g(x)|^{2p} \leq | \tilde f(x)|^{2q}$ for every $x \in U.$ Since $ \tilde f \in \hat\gta,$ we have that $ \tilde f^q \in \hat\gta.$ Let $f= \tilde f^q,$ and conclude that $|g(x)|^{2p} \leq | f(x)|^2$ for every $x \in U.$ We have arrived at the expression of the real radical used in \cite{abn} except for the fact that all computations here are on $U$ rather than on all of $\R^n.$ The rest of the proof given for Lemma 3.1 in \cite{abn} applies verbatim, so we refer the reader there for the rest of the argument. The hypothesis that $\ra(U) \subset \An(U)$ is used to approximate the coefficients in $\stb(U)$ of $f \in \hat\gta$ using elements of $\An(U)$ as every continuous function, hence every smooth one, can be approximated by a real analytic function. \qed

\smallskip\noindent {\bf Remarks:} 

\noindent (1)The conclusion that $  {\sqrt[\R]{\gta}}  =\sqrt[\R]{\hat\gta} \cap \An(U)$ would hold even if the real radical were computed germwise as in \ref{realrad} since the open set $U$ is chosen when the definition of the real radical is applied to $g$ and not modified subsequently. Therefore, Lemma~\ref{ra approx} can be used both on the presheaf and on the sheaf itself. 

\noindent (2) Note that by Proposition~\ref{allpropsthere} part (i), $\ra (U) \subset C^\R_M (U)$ and $\ra(U) \subset C_M (U),$ so we can let $\An(U) = C^\R_M (U)$ or $\An(U) = C_M (U).$

\noindent (3) Clearly, we can replace $U$ by any open neighborhood of $b \Omega \subset \C^n$ and apply this lemma to presheaves or sheaves defined on $b \Omega$ as in the Kohn algorithm.

\medskip\noindent In light of the previous three remarks, the significance of Lemma \ref{ra approx} is that it shows the operation $\wedge$ of tensoring with $\E$ commutes with the operation of taking the real radical $\sqrt[\R]{ \:}$ as follows:

$$
\begin{CD}
\gta     @>\wedge>>  \hat \gta \\@VV\sqrt[\R]{ \:}V       @VV\sqrt[\R]{ \:} V\\ \sqrt[\R]{ \gta}     @> \wedge>> \sqrt[\R]{ \hat \gta}=\widehat {\sqrt[\R]{ \gta}}
\end{CD}
$$
In other words, a priori tensoring with $\E$ and then taking the real radical would seem to produce a larger ideal than first taking the real radical and then tensoring with $\E.$ Lemma \ref{ra approx} shows that is not the case, so the previous diagram commutes.

\smallskip A priori, we would expect the full scale Kohn algorithm starting with a defining function $r$ in a Denjoy-Carleman class to produce many more subelliptic multipliers in $I^q_k$ for each $k$ than what we can obtain by working over just $C_M$ as in $\tilde I^q_k.$ The next theorem shows that is not the case. At each step $k,$ $I^q_k$ is obtained by taking all the Denjoy-Carleman subelliptic multipliers from the modified Kohn algorithm and generating over the ring of smooth functions. Since the sheaves $I^q_k$ were already proven to be quasicoherent in \cite{andreeaqf}, we can then use the good algebraic geometric properties of $\tilde I^q_k$ but work in $\widehat{\tilde I^q_k}$ where quasicoherence holds.

\medskip
\newtheorem{dcsheafsmooth}[raapprox]{Theorem}
\begin{dcsheafsmooth}
\label{dcsheafsmooththm} The ideal sheaf $I^q_k=\widehat{\tilde I^q_k}=\tilde I^q_k \otimes_{\E} \E$ for all $k \geq 1.$
\end{dcsheafsmooth}

\smallskip\noindent {\bf Remark:} This result is new even for the real analytic case discussed by Kohn in section 6 of \cite{kohnacta}. The sequence $M$ defining the Denjoy-Carleman class can be taken to consist of just $\{1, 1, 1, \dots \},$ which constructs the ring of real analytic functions $\ra.$ We can take $\An(U)$ to be $\ra(U)$ in Lemma~\ref{ra approx} and the proof will then be identical to the one given below.

\medskip\noindent {\bf Proof of Theorem~\ref{dcsheafsmooththm}:} $$\tilde I^q_1(x_0) = \sqrt[\R]{(\, r,\,
\text{coeff}\{\partial r \wedge \dbar r \wedge (\partial \dbar
r)^{n-q}\}\, )_{C_{M, b \Omega} (x_0)}},$$ 
whereas $$I^q_1(x_0) = \sqrt[\R]{(\, r,\,
\text{coeff}\{\partial r \wedge \dbar r \wedge (\partial \dbar
r)^{n-q}\}\, )_{\smooth_{b \Omega} (x_0)}}$$
The generators $r$ and $\text{coeff}\{\partial r \wedge \dbar r \wedge (\partial \dbar
r)^{n-q}\}$ are global Denjoy-Carleman sections in $C_M (b \Omega).$ The fact that tensoring with $\E$ commutes with the operation of taking the real radical implies that $I^q_1=\tilde I^q_1 \otimes_{\E} \E,$ so we have proven the result for $k=1.$

Next, we seek to show that the following diagram commutes:

$$
\begin{CD}
\tilde I^q_k    @>\text{advance step}>>  \tilde I^q_{k+1} \\@VV\wedge V       @VV\wedge V\\ \widehat{\tilde I^q_k}     @> \text{advance step}>> \widehat{\tilde I^q_{k+1}}
\end{CD}
$$
In other words, we seek to prove that tensoring with $\E$ commutes with the operation of moving to the next step in the Kohn algorithm. Recall that $$\tilde I^q_{k+1} (x_0) = \sqrt[\R]{(\, \tilde I^q_k
(x_0),\, \tilde A^q_k (x_0)\, )_{C_{M, b \Omega} (x_0)}},$$ where $$\tilde A^q_k (x_0)= \text{coeff}\{\partial
f_1 \wedge \dots \wedge \partial f_j \wedge \partial r \wedge
\dbar r \wedge (\partial \dbar r)^{n-q-j}\}$$ for $f_1, \dots, f_j
\in \tilde I^q_k (x_0)$ and $j \leq n-q.$ Since we have already proven in Lemma~\ref{ra approx} that tensoring with $\E$ commutes with the real radical, we only have to investigate what happens with $A^q_k (x_0)$ when we tensor with $\E.$ Let $f_i \in \widehat{\tilde I^q_k (x_0)}$ for $1 \leq i \leq j \leq n-q.$ Then $\displaystyle f_i = \sum_{l=1}^P \: h_l \, g_l$ for $g_l \in \tilde I^q_k (x_0)$ and $h_l \in \E(x_0)$ for all $1\leq l \leq P.$ Therefore, $\displaystyle  \partial f_i = \sum_{l=1}^P \: \partial h_l \, g_l+ h_l \partial g_l.$ When we compute $\partial f_1 \wedge \dots \wedge \partial f_j \wedge \partial r \wedge \dbar r \wedge (\partial \dbar r)^{n-q-j},$ the terms corresponding to $\partial h_l \, g_l$ where the smooth coefficient is being differentiated will give elements of $\widehat{\tilde I^q_k (x_0)}$ and only the terms corresponding to  $h_l \partial g_l$ will yield elements of $\widehat{\tilde A^q_k (x_0)}$ hence of $\widehat{\tilde I^q_{k+1} (x_0)}.$ We conclude that advancing the step in the Kohn algorithm does indeed commute with tensoring with $\E$ as needed. Inductively then, $I^q_k=\widehat{\tilde I^q_k}=\tilde I^q_k \otimes_{\E} \E$ for all $k \geq 1.$ \qed

\medskip\noindent We derive as an immediate corollary:

\smallskip
\newtheorem{dcsheafqc}[raapprox]{Corollary}
\begin{dcsheafqc}
\label{dcsheafqccor} The ideal sheaf $\widehat{\tilde I^q_k}$  is quasi-flasque hence quasicoherent for all $k \geq 1.$ 
\end{dcsheafqc}

\medskip\noindent {\bf Proof:} Since $I^q_k$ is quasi-flasque as proven in \cite{andreeaqf} and $I^q_k= \widehat{\tilde I^q_k}$ by Theorem~\ref{dcsheafsmooththm}, clearly $\widehat{\tilde I^q_k}$  is quasi-flasque hence quasicoherent for all $k \geq 1.$ \qed

\medskip As a consequence of Corollary~\ref{dcsheafqccor}, we can take advantage of the good algebraic geometric properties of $\tilde I^q_k$ as detailed in Section~\ref{DCdescription} but work with the quasicoherent sheaves $\widehat{\tilde I^q_k}.$

We now recall Theorem 3 proven by Eric Bedford and John Erik Forn\ae ss in 1981 in \cite{bf}. This theorem gives a generalization to $\smooth$ boundaries of the Diederich-Forn\ae ss theorem in \cite{df} that provides the
crucial geometrical step connecting the failure of the Kohn algorithm to generate the whole ring with
the existence of a complex variety in the boundary of the domain:

\medskip
\newtheorem{dfsmooth}[raapprox]{Theorem}
\begin{dfsmooth}
Let $\Omega \subset \C^n$ be a pseudoconvex domain with smooth boundary, and let $M \subset b \Omega$ be a smooth submanifold. If
$M$ has holomorphic dimension $q$ with respect to $b \Omega$ at some point $p \in M,$ then there exists a germ of a complex $q$-dimensional manifold with $V \subset b \Omega.$ Further, if $V$ cannot be chosen so that $V \cap M \neq \emptyset,$ then there is a manifold $V' \subset b \Omega$ with complex dimension $q+1.$
\label{dfsmooththm}
\end{dfsmooth}

\smallskip\noindent We are finally ready to tackle the proof of
the equivalence of types.

\medskip\noindent {\bf Proof of Theorem~\ref{maintheorem}:} As
explained in the introduction, the implication (iii) $\implies$
(ii) is the only one that needs to be proved. We will prove the
contrapositive statement that the failure of the Kohn algorithm to
terminate at the whole ring, negation of (ii), implies the
existence of a holomorphic variety in the boundary of the domain,
negation of finite D'Angelo type (iii). Let $x_0 \in b \Omega$ be
any point on the boundary of the domain. Consider the increasing chain of ideals given by the modified Kohn algorithm $\tilde I^q_1 (x_0) \subset \tilde I^q_2 (x_0) \subset \cdots$ By definition, for each $k \geq 1,$ $\tilde I^q_k(x_0) =  \sqrt[\R]{\tilde I^q_k(x_0)},$ so by the remark following the proof of Theorem~\ref{radaccprop} at the end of Section~\ref{DCdescription}, we conclude there exists some natural number $k$ such that $\tilde I^q_j(x_0) =\tilde I^q_k(x_0)$ for all $j \geq k.$ By
Proposition~\ref{finitetypeprop}, there exists a subideal of
finite type $\J_k \subset \tilde I_k^q (x_0)$ such that $\V(\J_k)$ can be taken as a representative of 
$\V(\tilde I^q_k(x_0)).$ We apply the {\L}ojasiewicz
Nullstellensatz, Theorem~\ref{nullstellensatz}, to $\J_k$ to
conclude that $\sqrt[\R]{\J_k}=\I(\V(\J_k))=\I(\V( \tilde I^q_k(x_0))).$
Since $\J_k \subset \tilde I^q_k(x_0) =  \sqrt[\R]{\tilde I^q_k(x_0)},$ the Nullstellensatz holds for $\tilde I^q_k(x_0)$
$$ \tilde I^q_k(x_0)=\I(\V( \tilde I^q_k(x_0))).$$ We needed to take this circuitous route because the ring $C_{M, b \Omega} (x_0)$ is not known to be Noetherian and the {\L}ojasiewicz
Nullstellensatz, Theorem~\ref{nullstellensatz}, only applies to finitely generated ideals. By assumption, the modified Kohn algorithm does not terminate, so $\tilde I^q_k(x_0)\neq C_{M, b \Omega} (x_0).$ Therefore, $\V(\tilde I^q_k(x_0))\neq \emptyset.$ Furthermore, $\V(\tilde I^q_k(x_0))$ cannot consist of just $x_0$ because otherwise, all functions $z_j$ would be multipliers and the algorithm would finish at the next step. As $\V(\tilde I^q_j(x_0)) = \V(\tilde I^q_k(x_0))$ for all $j \geq k,$
Proposition~\ref{stilltherethm}, Lemma~\ref{zariskireg}, and Proposition~\ref{regularptsprop} together imply that $ \V(\tilde I^q_k(x_0))$ has holomorphic dimension at least $q$ at $x_0$ itself and at a sequence of regular points that converges to $x_0.$ Let $f_1, \dots, f_p$ be the generators of $\J_k,$ which are chosen as in the proof of Proposition~\ref{regularptsprop} so that the support of the sheaf of finite type given by $(f_1, \dots, f_p)$ has the reduced structure. For each $i,$ $1 \leq i \leq p,$ there exists $U_i \ni x_0$ open in the topology of $b \Omega$ such that $f_i \in C_M (U_i).$ Set $\tilde U = U_1 \, \cap \dots \cap \, U_p.$ Clearly, $\tilde U$ is open, and from the proof of Proposition~\ref{regularptsprop}, $ \V(\tilde I^q_k(x_0)) \cap \tilde U$ contains a sequence of regular points converging to $x_0.$ As mentioned in the introduction, D'Angelo type is an open condition, so we shrink $\tilde U$ if necessary to also ensure that $b \Omega$ has finite D'Angelo type at every point of $\tilde U.$ We also shrink $\tilde U$ to make sure that the quasicoherence guaranteed by Corollary~\ref{dcsheafqccor} holds as well, i.e. that the Denjoy-Carleman sections generating $\widehat{\tilde I^q_k}(x_0)$ also generate $\widehat{\tilde I^q_k} (x)$ for every $x \in \tilde U.$ Next, we choose one of the regular points  $y \in \tilde U$ guaranteed to exist by Proposition~\ref{regularptsprop}, which has a neighborhood $U_y \subset \tilde U$ satisfying that $ \V(\tilde I^q_k(x_0)) \cap U_y$ is a manifold. We have already shown the holomorphic dimension of $ \V(\tilde I^q_k(x_0)) \cap U_y$ is at least $q$ at $y,$ but Lemma~\ref{zariskireg} and Proposition~\ref{stilltherethm} together yield that the same is true on all of $ \V(\tilde I^q_k(x_0)) \cap U_y$ up to potentially shrinking $U_y.$ Then Theorem~\ref{dfsmooththm} implies that there exists a complex manifold $\W \subset U_y$ such that $\dim_\C \W
\geq q.$ This contradicts condition (iii), finite D'Angelo type at every point of $\W.$ As the neighborhood $\tilde U \supset U_y \supset \W$ was chosen such that $b \Omega$ has finite D'Angelo type everywhere in $\tilde U,$ this is the needed contradiction.
\qed

\bigskip We now outline what crucial properties were required for Kohn's indirect argument as in section 6 of \cite{kohnacta} to go through in the Denjoy-Carleman case:
\begin{enumerate}[(i)]
\item Proving stabilization of the chain of ideals of multipliers in Kohn's algorithm was delicate even in the Denjoy-Carleman case. It is clear that it would be very difficult to do in the absence of at least topological Noetherianity, hence unlikely to hold in general in the smooth case given that the ring of smooth germs is not Noetherian and no topological Noetherianity is known.
\item  \L ojasiewicz inequalities are required for the real radical as defined by Kohn to yield a Nullstellensatz.
\item A Nullstellensatz is required to obtain equality in Lemma~\ref{zariski}, which is crucial in order to obtain a variety of holomorphic dimension $q$ in case the algorithm does not advance.
\item The Bedford-Forn\ae ss Theorem requires a smooth point, so smooth points must be dense, a property that is known to hold for \L ojasiewicz ideals, finitely generated ideals in the ring of smooth functions whose generators satisfy a \L ojasiewicz inequality. Perhaps the density of smooth points is true more generally, but it requires very good behavior as the variety is cut by hypersurfaces, something that the presence of a \L ojasiewicz inequality guarantees by default.
\item The functions cutting out the manifold at a smooth point used in the application of the Bedford-Forn\ae ss Theorem must be subelliptic multipliers, which once again requires a Nullstellensatz.
\item At least, quasicoherence holds for the sheaves of ideals of multipliers even in the smooth case as proven in \cite{andreeaqf}.
\end{enumerate}

\section*{Acknowledgements}
The author would like to thank Edward Bierstone and Pierre Milman for
introducing her to their work on the Denjoy-Carleman quasianalytic
classes and related literature as well as for very kindly taking
the time to answer numerous questions on the properties of these
classes. Additionally, the author would like to thank Vasile Brinzanescu and Fabrizio Broglia for a number of astute algebraic observations that led to important clarifications in the current manuscript.

\bibliographystyle{plain}
\bibliography{DCTypeEquiv}

\begin{thebibliography}{10}

\bibitem{nonwpt}
Francesca Acquistapace, Fabrizio Broglia, Michail Bronshtein, Andreea Nicoara,
  and Nahum Zobin.
\newblock Failure of the {W}eierstrass preparation theorem in quasi-analytic
  {D}enjoy-{C}arleman rings.
\newblock {\em Adv. Math.}, 258:397--413, 2014.

\bibitem{abn}
Francesca Acquistapace, Fabrizio Broglia, and Andreea~C. Nicoara.
\newblock A global {N}ullstellensatz for ideals of {D}enjoy-{C}arleman
  functions.
\newblock {\em Proc. Amer. Math. Soc.}, 144(5):2067--2071, 2016.

\bibitem{bf}
Eric Bedford and J.~E. Forn{\ae}ss.
\newblock Complex manifolds in pseudoconvex boundaries.
\newblock {\em Duke Math. J.}, 48(1):279--288, 1981.

\bibitem{bm}
Edward Bierstone and Pierre~D. Milman.
\newblock Resolution of singularities in {D}enjoy-{C}arleman classes.
\newblock {\em Selecta Math. (N.S.)}, 10(1):1--28, 2004.

\bibitem{bmv}
Edward Bierstone, Pierre~D. Milman, and Guillaume Valette.
\newblock Arc-quasianalytic functions.
\newblock {\em Proc. Amer. Math. Soc.}, 143(9):3915--3925, 2015.

\bibitem{boggess}
Albert Boggess.
\newblock {\em C{R} manifolds and the tangential {C}auchy-{R}iemann complex}.
\newblock Studies in Advanced Mathematics. CRC Press, Boca Raton, FL, 1991.

\bibitem{borel1}
\'{E}mile Borel.
\newblock Sur la g\'{e}n\'{e}ralisation du prolongement analytique.
\newblock {\em C. R. Acad. Sci. Paris}, 130:1115--1118, 1900.

\bibitem{borel2}
\'{E}mile Borel.
\newblock Sur les s\'{e}ries de polyn\^{o}mes et de fractions rationnelles.
\newblock {\em Acta Math.}, 24:309--387, 1901.

\bibitem{bazilandreea2019}
Vasile Brinzanescu and Andreea~C. Nicoara.
\newblock Relating {C}atlin and {D}'{A}ngelo {$q$}-types.
\newblock {\em J. Math. Anal. Appl.}, 527(1):Paper No. 127349, 13, 2023.

\bibitem{catlinnec}
David Catlin.
\newblock Necessary conditions for subellipticity of the {$\bar \partial
  $}-{N}eumann problem.
\newblock {\em Ann. of Math. (2)}, 117(1):147--171, 1983.

\bibitem{catlinbdry}
David Catlin.
\newblock Boundary invariants of pseudoconvex domains.
\newblock {\em Ann. of Math. (2)}, 120(3):529--586, 1984.

\bibitem{catlinsubell}
David Catlin.
\newblock Subelliptic estimates for the {$\overline\partial$}-{N}eumann problem
  on pseudoconvex domains.
\newblock {\em Ann. of Math. (2)}, 126(1):131--191, 1987.

\bibitem{catlindangelononeff}
David~W. Catlin and John~P. D'Angelo.
\newblock Subelliptic estimates.
\newblock In {\em Complex analysis}, Trends Math., pages 75--94.
  Birkh\"{a}user/Springer Basel AG, Basel, 2010.

\bibitem{childress}
C.~L. Childress.
\newblock Weierstrass division in quasianalytic local rings.
\newblock {\em Canad. J. Math.}, 28(5):938--953, 1976.

\bibitem{LeviCore2}
Gian~Maria Dall'Ara and Samuele Mongodi.
\newblock Remarks on the {L}evi core.
\newblock {\em Ann. Mat. Pura Appl. (4)}, 203(5):1997--2012, 2024.

\bibitem{opendangelo}
John~P. D'Angelo.
\newblock Real hypersurfaces, orders of contact, and applications.
\newblock {\em Ann. of Math. (2)}, 115(3):615--637, 1982.

\bibitem{dangelo95}
John~P. D'Angelo.
\newblock Finite type conditions and subelliptic estimates.
\newblock In {\em Modern methods in complex analysis ({P}rinceton, {NJ},
  1992)}, volume 137 of {\em Ann. of Math. Stud.}, pages 63--78. Princeton
  Univ. Press, Princeton, NJ, 1995.

\bibitem{df}
Klas Diederich and John~E. Fornaess.
\newblock Pseudoconvex domains with real-analytic boundary.
\newblock {\em Ann. Math. (2)}, 107(2):371--384, 1978.

\bibitem{hormest}
Lars H{\"o}rmander.
\newblock {$L\sp{2}$} estimates and existence theorems for the {$\bar \partial
  $}\ operator.
\newblock {\em Acta Math.}, 113:89--152, 1965.

\bibitem{kaplansky}
Irving Kaplansky.
\newblock {\em Commutative rings}.
\newblock The University of Chicago Press, Chicago, Ill.-London, revised
  edition, 1974.

\bibitem{kimzaitsev}
Sung-Yeon Kim and Dmitri Zaitsev.
\newblock Jet vanishing orders and effectivity of {K}ohn's algorithm in
  dimension 3.
\newblock {\em Asian J. Math.}, 22(3):545--568, 2018.

\bibitem{kimzaitsev2}
Sung-Yeon Kim and Dmitri Zaitsev.
\newblock Triangular resolutions and effectiveness for holomorphic subelliptic
  multipliers.
\newblock {\em Adv. Math.}, 387:Paper No. 107803, 36, 2021.

\bibitem{dbneumann1}
J.~J. Kohn.
\newblock Harmonic integrals on strongly pseudo-convex manifolds. {I}.
\newblock {\em Ann. of Math. (2)}, 78:112--148, 1963.

\bibitem{dbneumann2}
J.~J. Kohn.
\newblock Harmonic integrals on strongly pseudo-convex manifolds. {II}.
\newblock {\em Ann. of Math. (2)}, 79:450--472, 1964.

\bibitem{kohnweakpsc}
J.~J. Kohn.
\newblock Global regularity for {$\bar \partial $} on weakly pseudo-convex
  manifolds.
\newblock {\em Trans. Amer. Math. Soc.}, 181:273--292, 1973.

\bibitem{kohnacta}
J.~J. Kohn.
\newblock Subellipticity of the {$\bar \partial $}-{N}eumann problem on
  pseudo-convex domains: sufficient conditions.
\newblock {\em Acta Math.}, 142(1-2):79--122, 1979.

\bibitem{malgrange}
B.~Malgrange.
\newblock {\em Ideals of differentiable functions}.
\newblock Tata Institute of Fundamental Research Studies in Mathematics, No. 3.
  Tata Institute of Fundamental Research, Bombay, 1967.

\bibitem{andreeaqf}
Andreea~C. Nicoara.
\newblock Coherence and other properties of sheaves in the {K}ohn algorithm.
\newblock {\em Internat. J. Math.}, 25(8):1450077, 13, 2014.

\bibitem{racase}
Andreea~C. Nicoara.
\newblock Direct proof of termination of the {K}ohn algorithm in the
  real-analytic case.
\newblock {\em Pure Appl. Math. Q.}, 18(2):719--761, 2022.

\bibitem{op}
Jack Ohm and R.~L. Pendleton.
\newblock Rings with noetherian spectrum.
\newblock {\em Duke Math. J.}, 35:631--639, 1968.

\bibitem{pieroni}
Federica Pieroni.
\newblock On the real algebra of {D}enjoy-{C}arleman classes.
\newblock {\em Selecta Math. (N.S.)}, 13(2):321--351, 2007.

\bibitem{rsw}
J.-P. Rolin, P.~Speissegger, and A.~J. Wilkie.
\newblock Quasianalytic {D}enjoy-{C}arleman classes and o-minimality.
\newblock {\em J. Amer. Math. Soc.}, 16(4):751--777 (electronic), 2003.

\bibitem{siu2010}
Yum-Tong Siu.
\newblock Effective termination of {K}ohn's algorithm for subelliptic
  multipliers.
\newblock {\em Pure Appl. Math. Q.}, 6(4, Special Issue: In honor of Joseph J.
  Kohn. Part 2):1169--1241, 2010.

\bibitem{siu2017}
Yum-Tong Siu.
\newblock New procedure to generate multipliers in complex {N}eumann problem
  and effective {K}ohn algorithm.
\newblock {\em Sci. China Math.}, 60(6):1101--1128, 2017.

\bibitem{thilliez}
Vincent Thilliez.
\newblock On quasianalytic local rings.
\newblock {\em Expo. Math.}, 26(1):1--23, 2008.

\bibitem{thom}
Ren{\'e} Thom.
\newblock On some ideals of differentiable functions.
\newblock {\em J. Math. Soc. Japan}, 19:255--259, 1967.

\bibitem{tougeronqf}
Jean-Claude Tougeron.
\newblock Faisceaux diff\'erentiables quasi-flasques.
\newblock {\em C. R. Acad. Sci. Paris}, 260:2971--2973, 1965.

\bibitem{tougeron}
Jean-Claude Tougeron.
\newblock {\em Id\'eaux de fonctions diff\'erentiables}.
\newblock Springer-Verlag, Berlin, 1972.
\newblock Ergebnisse der Mathematik und ihrer Grenzgebiete, Band 71.

\end{thebibliography}

\end{document}